\def\date{March 14, 2007} 

\input amssym.def
\input amssym.tex

\def\item#1{\vskip1.3pt\hang\textindent {\rm #1}}


\newskip\litemindent
\litemindent=0.7cm  
\def\Litem#1#2{\par\noindent\hangindent#1\litemindent
\hbox to #1\litemindent{\hfill\hbox to \litemindent
{\ninerm #2 \hfill}}\ignorespaces}
\def\litem{\Litem1}

\tolerance=300
\pretolerance=200
\hfuzz=1pt
\vfuzz=1pt

\hoffset=0in
\voffset=0.5in

\hsize=5.8 true in 
\vsize=9.2 true in
\parindent=25pt
\mathsurround=1pt
\parskip=1pt plus .25pt minus .25pt
\normallineskiplimit=.99pt

\countdef\revised=100
\mathchardef\emptyset="001F 
\chardef\ss="19
\def\3{\ss}
\def\anf{$\lower1.2ex\hbox{"}$}
\def\frac#1#2{{#1 \over #2}}
\def\>{>\!\!>}
\def\<{<\!\!<}

\def\into{\hookrightarrow}
 
\def\ssssarr{\hbox to 15pt{\rightarrowfill}}
\def\sssarr{\hbox to 20pt{\rightarrowfill}}
\def\ssarr{\hbox to 30pt{\rightarrowfill}}
\def\sarr{\hbox to 40pt{\rightarrowfill}}
\def\arr{\hbox to 60pt{\rightarrowfill}}
\def\larr{\hbox to 60pt{\leftarrowfill}}
\def\Arr{\hbox to 80pt{\rightarrowfill}}

\def\sssmapright#1{\smash{\mathop{\sssarr}\limits^{#1}}}

\def\smapright#1{\smash{\mathop{\sarr}\limits^{#1}}}

\def\prolim{{\displaystyle \lim_{\longleftarrow}}\ }

\def\Ad{\mathop{\rm Ad}\nolimits}

\def\Aut{\mathop{\rm Aut}\nolimits}

\def\Evol{\mathop{\rm Evol}\nolimits}
\def\evol{\mathop{\rm evol}\nolimits}
\def\ev{\mathop{\rm ev}\nolimits}

\def\GL{\mathop{\rm GL}\nolimits}

\def\Hom{\mathop{\rm Hom}\nolimits}%
\def\id{\mathop{\rm id}\nolimits} 
\def\im{\mathop{\rm im}\nolimits}

\def\Im{\mathop{\rm Im}\nolimits}
\def\inf{\mathop{\rm inf}\nolimits}

\def\PGL{\mathop{\rm PGL}\nolimits}
\def\PSL{\mathop{\rm PSL}\nolimits}

\def\per{\mathop{\rm per}\nolimits}
\def\rad{\mathop{\rm rad}\nolimits}

\def\Re{\mathop{\rm Re}\nolimits}

\def\SL{\mathop{\rm SL}\nolimits}

\def\Spec{\mathop{\rm Spec}\nolimits}





\def\0{{\bf 0}}
\def\1{{\bf 1}}

\def\g{{\frak g}}
\def\gl{{\frak {gl}}}

\def\k{{\frak k}}

\def\sL{{\frak {sl}}}

\def\L{\mathop{\bf L{}}\nolimits}

\def\C{{{\Bbb C}{\mskip+1mu}}} 
\def\K{{{\Bbb K}{\mskip+2mu}}} 

\def\R{{\Bbb R}} 
\def\Z{{\Bbb Z}} 
\def\N{{\Bbb N}}

\def\K{{\Bbb K}}

\def\SS{{\Bbb S}} 
\def\T{{\Bbb T}} 

\def\:{\colon}  
\def\.{{\cdot}}
\def\|{\Vert}
\def\bsk{\bigskip}

\def\giantskip{\vskip2\bigskipamount}
\def\gsk{\giantskip}

\def\msk{\medskip}

\def \res {\!\mid\!\!}

\def\bbr{\bigbreak}
\def\giantbreak{\par \ifdim\lastskip<2\bigskipamount \removelastskip
         \penalty-400 \giantskip\fi}

\def\nin{\noindent}
\def\cen{\centerline}
\def\pagebreak{\vskip 0pt plus 0.0001fil\break}
\def\linebreak{\break}

\def\hat{\widehat}

\def\derat#1{{d \over dt} \hbox{\vrule width0.5pt 
                height 5mm depth 3mm${{}\atop{{}\atop{\scriptstyle t=#1}}}$}}

\def\eps{\varepsilon}
\def\epsilon{\varepsilon}
\def\eset{\emptyset}

\def\nin{\noindent}
\def\oline{\overline}

\def\pder#1,#2,#3 { {\partial #1 \over \partial #2}(#3)}
\def\pde#1,#2 { {\partial #1 \over \partial #2}}
\def\phi{\varphi}


\def\subeq{\subseteq}
\def\supeq{\supseteq}

\def\Rarrow{\Rightarrow}
\def\Larrow{\Leftarrow}
\def\tilde{\widetilde}

\font\ninerm=cmr9
\font\eightrm=cmr8

\font\eightbf=cmbx8


\font\smc=cmcsc10
\font\bfone=cmbx10 scaled\magstep1 
\font\bftwo=cmbx10 scaled\magstep2 

\def\qed{{\unskip\nobreak\hfil\penalty50\hskip .001pt \hbox{}\nobreak\hfil
          \vrule height 1.2ex width 1.1ex depth -.1ex
           \parfillskip=0pt\finalhyphendemerits=0\medbreak}\rm}

\def\qeddis{\eqno{\vrule height 1.2ex width 1.1ex depth -.1ex} $$
                   \medbreak\rm}

\def\Lemma #1. {\bigbreak\vskip-\parskip\noindent{\bf Lemma #1.}\quad\it}

\def\Sublemma #1. {\bigbreak\vskip-\parskip\noindent{\bf Sublemma #1.}\quad\it}

\def\Proposition #1. {\bigbreak\vskip-\parskip\noindent{\bf Proposition #1.}
\quad\it}

\def\Corollary #1. {\bigbreak\vskip-\parskip\nin{\bf Corollary #1.}
\quad\it}

\def\Theorem #1. {\bigbreak\vskip-\parskip\noindent{\bf Theorem #1.}
\quad\it}

\def\Definition #1. {\rm\bigbreak\vskip-\parskip\noindent
{\bf Definition #1.}
\quad}

\def\Remark #1. {\rm\bigbreak\vskip-\parskip\noindent{\bf Remark #1.}\quad}

\def\Example #1. {\rm\bigbreak\vskip-\parskip\noindent{\bf Example #1.}\quad}
\def\Examples #1. {\rm\bigbreak\vskip-\parskip\noindent{\bf Examples #1.}\quad}

\def\Problems #1. {\bigbreak\vskip-\parskip\noindent{\bf Problems #1.}\quad}
\def\Problem #1. {\bigbreak\vskip-\parskip\noindent{\bf Problem #1.}\quad}
\def\Exercise #1. {\bigbreak\vskip-\parskip\noindent{\bf Exercise #1.}\quad}

\def\Conjecture #1. {\bigbreak\vskip-\parskip\noindent{\bf Conjecture #1.}\quad}

\def\Proof#1.{\rm\par\ifdim\lastskip<\bigskipamount\removelastskip\fi\smallskip
            \noindent {\bf Proof.}\quad}

\def\Axiom #1. {\bigbreak\vskip-\parskip\noindent{\bf Axiom #1.}\quad\it}

\def\Satz #1. {\bigbreak\vskip-\parskip\noindent{\bf Satz #1.}\quad\it}

\def\Korollar #1. {\bbr\vskip-\parskip\nin{\bf Korollar #1.} \quad\it}

\def\Folgerung #1. {\bbr\vskip-\parskip\nin{\bf Folgerung #1.} \quad\it}

\def\Folgerungen #1. {\bbr\vskip-\parskip\nin{\bf Folgerungen #1.} \quad\it}

\def\Bemerkung #1. {\rm\bigbreak\vskip-\parskip\noindent{\bf Bemerkung #1.}
\quad}

\def\Beispiel #1. {\rm\bigbreak\vskip-\parskip\noindent{\bf Beispiel #1.}\quad}
\def\Beispiele #1. {\rm\bigbreak\vskip-\parskip\noindent{\bf Beispiele #1.}\quad}
\def\Aufgabe #1. {\rm\bigbreak\vskip-\parskip\noindent{\bf Aufgabe #1.}\quad}
\def\Aufgaben #1. {\rm\bigbreak\vskip-\parskip\noindent{\bf Aufgabe #1.}\quad}

\def\Beweis#1. {\rm\par\ifdim\lastskip<\bigskipamount\removelastskip\fi
           \smallskip\noindent {\bf Beweis.}\quad}

\nopagenumbers

\def\date{\ifcase\month\or January\or February \or March\or April\or May
\or June\or July\or August\or September\or October\or November
\or December\fi\space\number\day, \number\year}

\def\title{Title ??}
\def\author{Author ??}

\def\thanks#1{\footnote*{\eightrm#1}}

\def\rightheadline{\hfil{\eightrm\title}\hfil\tenbf\folio}
\def\leftheadline{\tenbf\folio\hfil{\eightrm\author}\hfil}
\headline={\vbox{\line{\ifodd\pageno\rightheadline\else\leftheadline\fi}}}

\def\firstheadline{}
\def\firstfootline{\cen{\rm\folio}}

\def\seite #1 {\pageno #1
               \headline={\ifnum\pageno=#1 \firstheadline
               \else\ifodd\pageno\rightheadline\else\leftheadline\fi\fi}
               \footline={\ifnum\pageno=#1 \firstfootline\else{}\fi}}

\newdimen\dimenone
 \def\checkleftspace#1#2#3#4{
 \dimenone=\pagetotal
 \advance\dimenone by -\pageshrink   
 \ifdim\dimenone>\pagegoal          
   \else\dimenone=\pagetotal
        \advance\dimenone by \pagestretch
        \ifdim\dimenone<\pagegoal
          \dimenone=\pagetotal
          \advance\dimenone by#1         
          \setbox0=\vbox{#2\parskip=0pt                
                     \hyphenpenalty=10000
                     \rightskip=0pt plus 5em
                     \noindent#3 \vskip#4}    
        \advance\dimenone by\ht0
        \advance\dimenone by 3\baselineskip   
        \ifdim\dimenone>\pagegoal\vfill\eject\fi
          \else\eject\fi\fi}


\def\subheadline #1{\nin\bigbreak\vskip-\lastskip
      \checkleftspace{0.9cm}{\bf}{#1}{\medskipamount}
          \indent\vskip0.7cm\centerline{\bf #1}\medskip}
\def\subsection{\subheadline} 

\def\lsubheadline #1 #2{\bigbreak\vskip-\lastskip
      \checkleftspace{0.9cm}{\bf}{#1}{\bigskipamount}
         \vbox{\vskip0.7cm}\cen{\bf #1}\msk \cen{\bf #2}\bsk}

\def\sectionheadline #1{\bigbreak\vskip-\lastskip
      \checkleftspace{1.1cm}{\bf}{#1}{\bigskipamount}
         \vbox{\vskip1.1cm}\cen{\bfone #1}\bsk}
\def\section{\sectionheadline} 

\def\lsectionheadline #1 #2{\bigbreak\vskip-\lastskip
      \checkleftspace{1.1cm}{\bf}{#1}{\bigskipamount}
         \vbox{\vskip1.1cm}\cen{\bfone #1}\msk \cen{\bfone #2}\bsk}

\def\lchapterheadline #1 #2{\bigbreak\vskip-\lastskip\indent\vskip3cm
                       \cen{\bftwo #1} \msk \cen{\bftwo #2} \gsk}
\def\llsectionheadline #1 #2 #3{\bigbreak\vskip-\lastskip\indent\vskip1.8cm
\cen{\bfone #1} \msk \cen{\bfone #2} \msk \cen{\bfone #3} \nobreak\bsk\nobreak}


\newtoks\literat
\def\[#1 #2\par{\literat={#2\unskip.}%
\hbox{\vtop{\hsize=.15\hsize\nin [#1]\hfill}
\vtop{\hsize=.82\hsize\nin\the\literat}}\par
\vskip.3\baselineskip}

\def\references{
\sectionheadline{\bf References}
\frenchspacing

\entries\par}

\mathchardef\emptyset="001F 
\def\address{Author: \tt$\backslash$def$\backslash$address$\{$??$\}$}

\def\abstract #1{{\narrower\baselineskip=10pt{\noindent
\eightbf Abstract.\quad \eightrm #1 }
\bigskip}}

\def\firstpage{\nin
{\obeylines \parindent 0pt }
\vskip2cm
\centerline{\bfone\title}
\gsk
\centerline{\bf\author}
\vskip1.5cm \rm}

\def\addresstwo{}

\def\dlastpage{\par\vbox{\vskip1cm\nin
\line{
\vtop{\hsize=.5\hsize{\parindent=0pt\baselineskip=10pt\nin\address}}
\quad 
\vtop{\hsize=.42\hsize\nin{\parindent=0pt
\baselineskip=10pt\addresstwo}}
\hfill} }}

\def\Box #1 { \msk\par\nin 
\centerline{
\vbox{\offinterlineskip
\hrule
\hbox{\vrule\strut\hskip1ex\hfil{\smc#1}\hfill\hskip1ex}
\hrule}\vrule}\msk }

\def\adots{\mathinner{\mkern1mu\raise1pt\vbox{\kern7pt\hbox{.}}
                        \mkern2mu\raise4pt\hbox{.}
                        \mkern2mu\raise7pt\hbox{.}\mkern1mu}}


\pageno=1

\def\shtitle{Lie group structures on groups of maps} 
\def\title{
\vbox{\centerline{Lie group structures on groups of smooth and holomorphic maps} 
\centerline{on non-compact manifolds}}}
\def\author{}
\def\author{Karl-Hermann Neeb, Friedrich Wagemann}
\def\leftheadline{\tenbf\folio\hfil{\tt ncc8.tex}\hfil\eightrm\date}
\def\rightheadline{\tenbf\folio\hfil{\rm\shtitle}\hfil\eightrm\date}

\hyphenation{Gro-then-dieck} 
\def\MC{\mathop{\rm MC}\nolimits}

\def\address{Karl-Hermann Neeb 

Fachbereich Mathematik 

Technische Universit\"at Darmstadt 

Schlossgartenstr. 7 

64285 Darmstadt 

Germany

neeb@mathematik.tu-darmstadt.de} 

\def\addresstwo{Friedrich Wagemann 

     Laboratoire de Math\'ematiques 

Jean Leray

     Facult\'e des Sciences et Techniques

     Universit\'e de Nantes

     2, rue de la Houssini\`ere

     44322 Nantes cedex 3

     France 

wagemann@math.univ-nantes.fr}

\def\sst{\scriptstyle}
\firstpage 

\abstract{We study Lie group structures on groups of the form 
${\sst C^\infty(M,K)}$, where ${\sst M}$ 
is a non-compact smooth manifold and ${\sst K}$ 
is a, possibly infinite-dimensional, Lie group. 
First we prove that there is at most one 
Lie group structure with Lie algebra ${\sst C^\infty(M,\k)}$ for which the 
evaluation map is smooth. We then prove the existence of such a structure  
if the universal cover of ${\sst K}$ is diffeomorphic to a locally convex space 
and if the image of the left logarithmic derivative in 
${\sst \Omega^1(M,\k)}$ 
is a smooth submanifold, the latter being 
the case in particular if ${\sst M}$ is one-dimensional. We also obtain 
analogs of these results for the group ${\sst {\cal O}(M,K)}$ of holomorphic 
maps on a complex manifold with values in a complex Lie group~${\sst K}$. 
We further show that there exists a natural Lie group 
structure on ${\sst {\cal O}(M,K)}$ 
if ${\sst K}$ is Banach and ${\sst M}$ is a non-compact complex curve 
with finitely generated fundamental group. 
\hfill\linebreak 
AMS Classification: 22E65, 22E67, 22E15, 22E30  \hfill\linebreak 
}

\sectionheadline{Introduction} 

\nin If $M$ is a finite-dimensional manifold (possibly with boundary) 
and $K$ a Lie group (modeled on a locally convex space), 
then the
group $C^\infty(M,K)$ of smooth maps with values in $K$ has a 
natural group topology, called the {\it smooth compact open topology}. 
If, in addition, $M$ is a complex manifold without boundary and 
$K$ is a complex Lie group, then $C^\infty(M,K)$ contains the subgroup 
${\cal O}(M,K)$ of holomorphic maps, on which the smooth compact open 
topology simply 
coincides with the compact open topology if $M$ has no boundary. 
In this paper we discuss the 
question when the topological groups $C^\infty(M,K)$, resp. ${\cal O}(M,K)$, 
carry Lie group structures. 

If $K = E$ is a locally convex space, then $C^\infty(M,E)$ also is a 
locally convex space, hence a Lie group. It is also well-known that 
if $M$ is compact, then $C^\infty(M,K)$ carries a natural Lie group 
structure (cf.\ [Mi80], [Gl02a]; and [Wo06] for manifolds with boundary). 
In this case one obtains charts of $C^\infty(M,K)$ 
by composing with charts of $K$. This does no longer work for non-compact manifolds. 
A necessary condition for a topological group $G$ to possess a compatible Lie group
structure is that it is locally contractible, which implies in
particular that the arc-component $G_a$ of the identity is open. 
In some cases we shall prove that the topological group 
$C^\infty(M,K)$ is not a Lie group by 
showing that the latter condition fails. 

In general, we cannot expect the group 
$C^\infty(M,K)$ to carry a Lie group structure, but any 
``reasonable'' Lie group structure on this group should have the property that 
for any smooth manifold $N$, a map $f \: N \to C^\infty(M,K)$ is smooth if and only if 
the corresponding map 
$$ f^\wedge \: N \times M \to K, \quad f^\wedge(n,m) := f(n)(m) $$
is smooth. We thus start our investigation in Section~I 
with a characterization of 
smooth maps $N \times M \to K$ in terms of data associated to the group 
$G := C^\infty(M,K)$. This leads to the main result of Section I, that 
for any regular Lie group $K$ the group $G$ carries at 
most one regular Lie group structure with Lie algebra $\g = C^\infty(M,\k)$ 
for which all evaluation maps $\ev_m \: C^\infty(M,K) \to K$ are smooth 
with $\L(\ev_m) = \ev_m$. We call such a Lie group structure {\it compatible 
with evaluations}. 

>From what we have said above, it easily follows that 
such Lie group structures exist if 
$M$ is compact (Theorem~I.3) or if the universal covering group 
$\tilde K$ of $K$ is diffeomorphic to a locally convex space 
(Theorem~IV.2), which is the case if $K$ is regular 
abelian or finite-dimensional solvable. If $\tilde K$ is diffeomorphic to 
a locally convex space $E$, 
the Lie group  $C^\infty(M,\tilde K)$ is diffeomorphic to the 
locally convex space $C^\infty(M,E)$ and the underlying topology coincides
with the smooth compact open topology. If $K$ is not simply connected, 
the Lie topology might be finer 
than the smooth compact open topology, but both coincide on the 
arc-component of the identity, which need not be open in the smooth 
compact open topology. In Section IV we take a closer look at this subtle 
situation and show that if, f.i., $K$ is finite-dimensional solvable, 
then both topologies coincide if and only if the group $H^1(M,\Z)$ is finitely 
generated (Remark~IV.13). 
We also derive analogous results for the group 
${\cal O}(M,K)$ of holomorphic maps on a complex manifold $M$ with values 
in a complex Lie group $K$. 

Clearly, the condition that $\tilde K$ is diffeomorphic to a vector space is 
quite restrictive. 
To find weaker sufficient conditions for the existence of Lie group 
structures, we study in Section II the left logarithmic derivative 
$$ \delta \: C^\infty(M,K) \to \MC(M,\k) := \{ \alpha \in \Omega^1(M,\k) \: 
d\alpha + {\textstyle{1\over 2}}[\alpha,\alpha]= 0\}, \quad 
f \mapsto f^{-1}.df. $$
We show in Proposition~II.1 that for any regular Lie group $K$ 
with Lie algebra $\k$, any connected manifold $M$ and any $m_0 \in M$, it maps 
the subgroup 
$$C^\infty_*(M,K) := \{ f \in C^\infty(M,K) \: f(m_0) = \1\}$$ 
of based maps homeomorphically onto its image, which is characterized by the 
Fundamental Theorem (Theorem~I.5). 
If $\im(\delta)$ carries a natural manifold 
structure, we thus obtain a manifold structure on the group $C^\infty_*(M,K)$ and 
hence on $C^\infty(M,K) \cong K \ltimes C^\infty_*(M,K)$ a regular 
Lie group structure compatible with evaluations (Theorem~II.2). 
If $M$ is $1$-connected, this is the case if 
$K$ is abelian, $M$ is real and one-dimensional 
(cf.\ [KM97]) or for holomorphic maps on complex curves. 

If $M$ is not simply connected, the situation is more complicated. 
If $K$ is abelian, we always have a regular Lie group structure on 
$C^\infty(M,K)$ since $\tilde K$ is a vector space, but it is compatible 
with the smooth compact open topology only if $H^1(M,\Z)$ is finitely 
generated (Theorem~IV.8). If $M$ is real $1$-dimensional, then it either 
is compact or simply connected, but the situation becomes interesting 
if $M$ is a complex curve which is not simply connected. 
It turns out that in this situation one can show that $\delta({\cal O}(M,K))$ 
is a complex submanifold of the Fr\'echet space $\Omega^1_h(M,\k)$ of 
holomorphic $\k$-valued $1$-forms on $M$ whenever 
$K$ is a Banach--Lie group and $\pi_1(M)$ 
is finitely generated. The key tool in our argument is  
Gl\"ockner's Implicit Function Theorem ([Gl03]) for smooth maps 
on locally convex spaces with values in Banach spaces, which is used 
to take care of the period conditions. 

We collect some of our main results in the following two theorems: 
\Theorem 1. Let $K$ be a connected regular real Lie group and $M$ 
a real finite-dimensional connected manifold. 
Then the group $C^\infty(M,K)$ carries a Lie group structure compatible with 
evaluations if 
\litem{(1)} $\tilde K$ is diffeomorphic to a locally convex space. 
If, in addition, $\pi_1(M)$ is finitely generated, the Lie group structure 
is compatible with the smooth compact open topology {\rm(Theorem~IV.2)}. 
\litem{(2)} $\dim M = 1$ {\rm(Corollary~II.3)}. 
\litem{(3)} $M \cong \R^k \times C$, where $C$ is compact {\rm(Corollary~II.8)}. 
\qed
   
For complex groups and holomorphic maps we have: 

\Theorem 2. Let $K$ be a regular complex Lie group and $M$ 
a finite-dimensional connected complex manifold without boundary. 
Then the group ${\cal O}(M,K)$ carries a Lie group structure with 
Lie algebra ${\cal O}(M,\k)$ compatible with evaluations if
\litem{(1)} $\tilde K$ is diffeomorphic to a locally convex space. If,  
in addition, $\pi_1(M)$ is finitely generated, the Lie group structure 
is compatible with the compact open topology {\rm(Theorem~IV.3)}. 
\litem{(2)} $\dim_\C M = 1$, $\pi_1(M)$ is finitely generated 
and $K$ is a Banach--Lie group {\rm(Theorem~III.12)}. 
\qed

Actually (2) in the preceding theorem was the original source of 
motivation for this work. It provides in particular a Lie theoretic environment 
for Lie groups associated to Krichever-Novikov Lie algebras which form an 
interesting generalization of affine Kac--Moody algebras (cf.\ [Sch03]).

If $M$ is a $\sigma$-compact finite-dimensional manifold 
and $M= \bigcup_n M_n$ is an exhaustion of $M$ by compact submanifolds 
$M_n$ with boundary, then the group 
$C^\infty(M,K)$ can be identified with the projective limit 
$\prolim C^\infty(M_n,K)$, where the connecting maps are given by 
restriction. Since each group $C^\infty(M_n,K)$ carries a natural Lie group 
structure, the topological group $C^\infty(M,K)$ is 
a projective limit of Lie groups. From this point of view, the present 
paper deals with Lie group structures on certain projective limits 
of infinite-dimensional Lie groups.  
For projective limits 
of finite-dimensional Lie groups, the problem to characterize the Lie groups 
among these groups has been solved completely in [HoNe06]. 

The paper is structured as follows. Section I contains generalities on 
smooth maps with values in Lie groups and the aforementioned uniqueness 
result on Lie group structures with smooth evaluation map 
(Corollary~I.10). In Section II we exploit the method to obtain 
Lie group structures on $C^\infty(M,K)$ by submanifold structures on 
$\im(\delta)$ and in Section III we transfer this method to groups 
of holomorphic maps. Target groups $K$ whose universal cover is diffeomorphic 
to a locally convex space are discussed in Section~IV. 
We conclude with a short section on strange properties of the exponential 
map of groups of holomorphic maps on non-compact manifolds and an 
appendix with some technical tools necessary to deal with 
manifolds of smooth and holomorphic maps. 

{\bf Acknowledgement:} We are most grateful to H.\ Gl\"ockner for 
many useful comments on earlier versions of the manuscript.

\subheadline{Preliminaries} 

Let $X$ and $Y$ be locally convex topological vector spaces, $U
\subeq X$ open and $f \: U \to Y$ a map. Then the {\it derivative
  of $f$ at $x$ in the direction of $h$} is defined as 
$$ df(x)(h) := \lim_{t \to 0} {1 \over t} \big( f(x + t h) - f(x)\big)
$$
whenever the limit exists. The function $f$ is called {\it differentiable at
  $x$} if $df(x)(h)$ exists for all $h \in X$. It is called {\it
  continuously differentiable or $C^1$} if it is continuous and differentiable at all
points of $U$ and 
$$ df \: U \times X \to Y, \quad (x,h) \mapsto df(x)(h) $$
is a continuous map. It is called a $C^n$-map if it is $C^1$ and $df$ is a
$C^{n-1}$-map, and $C^\infty$ (or {\it smooth}) if it is $C^n$ for all $n \in \N$. 
This is the notion of differentiability used in [Mil84], [Ha82] and
[Gl02b], where the latter reference deals with the modifications
needed for incomplete spaces. If $X$ and $Y$ are complex, 
$f$ is called {\it holomorphic} if it is smooth and its differentials $df(x)$ are 
complex linear. If $Y$ is Mackey complete, it suffices that $f$ is $C^1$. 

Since we have a chain rule for $C^1$-maps between locally convex 
spaces, we can define smooth manifolds $M$ as in
the finite-dimensional case. A chart $(\phi,U)$ with respect to a given 
manifold structure on $M$ is an open set $U\subset M$ together with a 
homeomorphism $\phi$ onto an open set of the model space. 
An atlas for the tangent bundle $TM$ is obtained
directly from an atlas of $M$, but we do not consider the cotangent
bundle as a manifold because this requires to choose a topology on the
dual spaces, for which there are many possibilities. 
Nonetheless, there is a natural concept of a smooth $k$-form on $M$. 
If $E$ is a locally convex space, then an {\it $E$-valued $k$-form} $\omega$ on $M$ is a function 
$\omega$ which associates to each $p \in M$ a $k$-linear 
alternating map $T_p(M)^k \to E$ such that in local coordinates the
map 
$(p,v_1, \ldots, v_k) \mapsto \omega(p)(v_1, \ldots, v_k)$
is smooth. We write $\Omega^k(M,E)$ for the space of smooth $k$-forms
on $M$ with values in $E$. The differentials 
$$ d \: \Omega^k(M,E) \to \Omega^{k+1}(M,E) $$
are defined by the same formula as in the finite-dimensional case 
(cf.\ [Beg87]).

If $M$ is a smooth manifold modeled on the locally convex space $E$,  
a subset $N \subeq M$ is called a {\it submanifold} 
of $M$ if there exists a closed subspace 
$F \subeq E$ and for each $n \in N$ there exists an $E$-chart $(\phi,U)$ of 
$M$ with $n \in U$ and $\phi(U \cap N) = \phi(U) \cap F$. 
The submanifold $N$ is called a {\it split submanifold} if, in addition, 
there exists a subspace $G \subeq E$ for which the addition map 
$F \times G \to E, (f,g) \mapsto f + g$ is a topological isomorphism.

$M$ is a smooth manifold modeled on the locally convex space $E$ 
with boundary $\partial M$ in case the 
$m\in \partial M$ have smooth charts $(\phi,U)$ to open neighborhoods $\phi(U)$
of the boundary of a half space of $E$. $M$ is said to have corners in case 
corner points have smooth charts to open neighborhoods of the vertex of a 
quadrant in $E$. Boundaries and the set of corners may be empty, and
$M$ reduces in this case to an ordinary manifold. For a complex manifold $M$,
the boundary is always supposed to be empty. 

A {\it Lie group} $G$ is a group equipped with a 
smooth manifold structure modeled on a locally convex space 
for which the group multiplication and the 
inversion are smooth maps. We write $\1 \in G$ for the identity element and 
$\lambda_g(x) = gx$, resp., $\rho_g(x) = xg$ for the left, resp.,
right multiplication on $G$. Then each $x \in T_\1(G)$ corresponds to
a unique left invariant vector field $x_l$ with 
$x_l(g) := d\lambda_g(\1).x, g \in G.$
The space of left invariant vector fields is closed under the Lie
bracket of vector fields, hence inherits a Lie algebra structure. In
this sense we obtain on $\g := T_\1(G)$ a continuous Lie bracket which
is uniquely determined by $[x,y]_l = [x_l, y_l]$ for $x,y \in \g$. 
The {\it Maurer--Cartan form} $\kappa_G \in \Omega^1(G,\g)$ is the unique 
left invariant $1$-form on $G$ with $\kappa_{G,\1} = \id_\g$, i.e., 
$\kappa_G(x_l) = x$ for each $x \in \g$. We write 
$q_G \: \tilde G_0 \to G_0$ for the universal covering map of the identity 
component $G_0$ of $G$ and identify 
the discrete central subgroup $\ker q_G$ of $\tilde G_0$ with $\pi_1(G) \cong \pi_1(G_0)$. 

In the following we always write $I = [0,1]$ for the unit interval in $\R$. 
A Lie group $G$ is called {\it regular} if for each 
$\xi \in C^\infty(I,\g)$, the initial value problem 
$$ \gamma(0) = \1, \quad \gamma'(t) = \gamma(t).\xi(t) = T(\lambda_{\gamma(t)})\xi(t) $$
has a solution $\gamma_\xi \in C^\infty(I,G)$, and the
evolution map 
$$ \evol_G \: C^\infty(I,\g) \to G, \quad \xi \mapsto \gamma_\xi(1) $$
is smooth (cf.\ [Mil84]). We then also write 
$$ \Evol_G \: C^\infty(I,\g) \to C^\infty(I,G), \quad 
\xi \mapsto \gamma_\xi, $$
and recall that this is a smooth map if $C^\infty(I,G)$ carries its natural 
Lie group structure (cf.\ Theorem~I.3 and Lemma~A.5 below). 
For a locally convex space $E$, the regularity of the Lie group $(E,+)$ is equivalent 
to the Mackey completeness of $E$, i.e., to the existence of integrals 
of smooth curves $\gamma \: I \to E$. We also recall that for each regular 
Lie group $G$, its Lie algebra $\g$ is regular and that all Banach--Lie groups 
are regular ([GN07]). The evolution map 
$\evol_G$ is supposed to be holomorphic for a complex regular Lie group $G$.

Throughout this paper, $K$ denotes a regular Lie group. 

\sectionheadline{I. Generalities on groups of smooth maps and regular Lie groups} 

In this section we introduce the natural group topology on the group 
$G = C^\infty(M,K)$ of smooth 
maps from a manifold $M$ with values in a Lie group $K$. 
We then describe some technical tools to deal 
with nonlinear maps between spaces of smooth maps and differential forms 
which leads to a characterization of maps $f \: N \to G$ for which the corresponding 
map $f^\wedge \: N \times M \to K$ is smooth (Proposition~I.8). This in turn is used 
to show that $G$ carries at most one regular Lie group structure 
compatible with evaluations (Corollary~I.10).

\Definition I.1. (Groups of differentiable maps as topological groups)  
(a) If $X$ and $Y$ are topological spaces, then the {\it compact open topology} 
on the space $C(X,Y)$ is defined as  the topology generated by the sets of the form 
$$ W(K,U) := \{ f \in C(X,Y) \: f(K) \subeq U \}, $$
where $K$ is a compact subset of $X$ and $U$ an open subset of $Y$. 
We write $C(X,Y)_c$ for the topological space obtained by endowing 
$C(X,Y)$ with the compact open topology. 

(b) If $K$ is a topological group and $X$ is Hausdorff, then $C(X,K)$ is a group with respect to the 
pointwise product. Then the compact open topology on $C(X,K)$ coincides 
with the topology of uniform convergence on compact subsets of $X$, for which 
the sets $W(C,U)$, $C \subeq X$ compact and $U \subeq K$ a $\1$-neighborhood in $K$,  
form a basis of $\1$-neighborhoods. In particular, 
$C(X,K)_c$ is a topological group. 

(c) In the following we topologize for two smooth manifolds $M$ (possibly with 
boundary) and $N$, the space $C^\infty(M,N)$ by the embedding 
$$ C^\infty(M,N) \into \prod_{k = 0}^\infty C(T^k(M),T^k(N))_c, \quad 
f \mapsto (T^k(f))_{k \in \N_0}, \leqno(1.1) $$
where the spaces $C(T^k(M),T^k(N))_c$ carry the compact open topology. 
The so obtained topology on $C^\infty(M,N)$ is called the {\it smooth compact 
open topology}. 

Now let $K$ be a Lie group with Lie algebra $\k$ and $r \in \N_0 \cup \{\infty\}$. 
The tangent map $T(m_K)$ of the multiplication map 
$m_K \: K \times K \to K$ defines a Lie group structure 
on the tangent bundle $TK$ (cf.\ [GN07]). 
Iterating this procedure, we obtain a Lie group structure on all higher tangent bundles $T^n K$. 
For each $n \in \N_0$, we thus obtain topological groups $C(T^n M, T^n K)_c$. 
We also observe that for two smooth maps $f_1, f_2 \: M \to K$,  
the functoriality of $T$ yields 
$$ T(f_1\cdot f_2) = T(m_K \circ (f_1 \times f_2)) 
=T(m_K) \circ (Tf_1 \times Tf_2) = Tf_1 \cdot Tf_2. $$
Therefore the inclusion map 
$C^\infty(M,K) \into \prod_{n = 0}^\infty  C(T^n M, T^n K)_c$ 
from (1.1) is a group homomorphism, so that the inverse image of the product 
topology from the right hand side is a 
group topology on $C^\infty(M,K)$, called the {\it smooth compact open topology}. 
It turns $C^\infty(M,K)$ into a topological group, 
even if $M$ and $K$ are infinite-dimensional. 

(d) In the following we topologize the space $\Omega^1(M,E)$ of 
$E$-valued $1$-forms on $M$ as a closed subspace of 
$C^\infty(TM,E)$. 
\qed

For later reference, we first collect some information on the case 
where $K = E$ is a locally convex space or where $M$ is compact. 

\Proposition I.2. Let $M$ be a finite-dimensional 
smooth manifold and $E$ a locally convex space. 
Then the following assertions hold: 
\litem{(1)} $C^\infty(M,E)$ is a locally convex space, hence a Lie group 
and the evaluation map of $C^\infty(M,E)$ is smooth. 
If $E$ is Mackey complete, then $C^\infty(M,E)$ is Mackey complete, 
hence a regular Lie group. 
\litem{(2)} If $M$ and $E$ are complex, then 
${\cal O}(M,E) \into C^\infty(M,E)$ is a closed subspace, 
and the  evaluation map 
$\ev \: {\cal O}(M,E) \times M \to E$
is holomorphic. 
If $E$ is Mackey complete, then ${\cal O}(M,E)$ is Mackey complete, 
hence a regular Lie group. If $M$ has no boundary, then 
the subspace topology on ${\cal O}(M,E)$ coincides with 
the compact open topology. 

\Proof. (1) All the spaces $C(T^k M, T^k E)_c$ are 
locally convex. Therefore the corresponding product topology is 
locally convex, and hence $C^\infty(M,E)$ is a locally convex space. 

The continuity 
of the evaluation map follows from the continuity of the evaluation map for the 
compact open topology because the topology on $C^\infty(M,E)$ is finer. 
Next we observe that directional derivatives exist and lead to a map 
$$ d\ev \: C^\infty(M,E)^2 \times T(M) \to E, \quad 
((f,\xi),v_m) \mapsto \xi(m) + T_m(f) v $$
whose continuity follows from the first step, applied to the evaluation map of 
$C^\infty(TM,E)$. Hence $\ev$ is $C^1$, and iteration of this argument yields 
smoothness. 

In view of Lemma A.3, we have 
$C^\infty(I, C^\infty(M,E)) \cong C^\infty(I \times M, E),$
and if $E$ is Mackey complete, then we have an integration map 
$$ C^\infty(I \times M, E) \to C^\infty(M,E), \quad \xi \mapsto 
\int_0^1 \xi(t,\cdot)\, dt $$
which implies that each smooth curve with values in 
$C^\infty(M,E)$ has a Riemann integral, i.e., that 
$C^\infty(M,E)$ is Mackey complete, which, in view of $\evol(\xi) = \int_0^1 \xi(t)\, dt$, 
is equivalent to it being a regular Lie group. 

(2) Since ${\cal O}(M,E)$ is a closed subspace of $C^\infty(M,E)$, 
${\cal O}(M,E) \times M$ is a closed submanifold of 
$C^\infty(M,E) \times M$, and the first part of the proof shows that 
the evaluation map is smooth on this space. Clearly, it is separately 
holomorphic in both arguments, hence its differential is complex linear, 
and the assertion follows. 

If $E$ is Mackey complete, then $C^\infty(M,E)$ is Mackey complete by (1), 
and the closed subspace ${\cal O}(M,E)$ inherits this property. 

If $M$ has no boundary, then the Cauchy Formula 
entails that on the space ${\cal O}(M,E)$ uniform convergence on compact 
subsets implies in any local chart uniform convergence 
of all partial derivatives on compact subsets. Hence the inclusion map 
$$ {\cal O}(M,E) \into C^\infty(M,E) $$
is continuous and therefore a topological 
embedding.\footnote{$^1$}{\eightrm Note that if ${\sst M}$ is a complex manifold with boundary, then we cannot 
expect that the topology ${\sst {\cal O}(M,\C)}$ inherits from ${\sst 
C^\infty(M,\C)}$ 
coincides with the compact open topology, as can be seen for the example 
${\sst M = \{ z \in \C \: |z| < 1\}}$. In this case the space 
${\sst {\cal O}(M,\C)}$ of holomorphic functions with smooth boundary values 
is not complete with respect to the compact open topology, but it is a closed 
subspace of ${\sst C^\infty(M,\C)}$ with respect to the smooth compact open topology.}
\qed

\Theorem I.3. Let $M$ be a smooth manifold 
and $K$ be a Lie group with Lie algebra $\k$. Then the following 
assertions hold: 
\litem{(1)} If $M$ is compact (possibly with corners or boundary), 
then $C^\infty(M,K)$ carries a Lie group structure for which any $\k$-chart 
$(\phi_K,U_K)$ of $K$ yields a $C^\infty(M,\k)$-chart $(\phi,U)$ with 
$$ U := \{ f \in C^\infty(M,K) \: f(M) \subeq U_K\}, \quad 
\phi(f) := \phi_K \circ f,  $$
and the evaluation map of $C^\infty(M,K)$ is smooth. 
The corresponding Lie algebra is \break $C^\infty(M,\k)$, and if $K$ is regular, 
then $C^\infty(M,K)$ is regular. 
\litem{(2)} If $M$ is compact and complex (possibly with boundary)  
and $K$ is a complex Lie group, 
then ${\cal O}(M,K)$, endowed with the smooth compact open topology, 
carries a Lie group structure for which any chart 
$(\phi_K,U_K)$ of $K$ yields a chart $(\phi,U)$ with 
$$ U := \{ f \in {\cal O}(M,K) \: f(M) \subeq U_K\}, \quad 
\phi(f) := \phi_K \circ f,  $$
and the evaluation map $\ev \: {\cal O}(M,K) \times M \to K$ 
is holomorphic. The corresponding Lie algebra is ${\cal O}(M,\k)$.

\Proof. (1) For the existence of the Lie group structure with the given charts we 
refer to [Gl02a] for the case without boundary which is also dealt 
with in  [Mi80], and to [Wo05] for the case of manifolds with corners, 
including in particular manifolds with boundary. 

The smoothness of the evaluation map follows on each domain 
$U$ as above from the openness of $C^\infty(M,\phi_K(U_K))$ in 
$C^\infty(M,\k)$ and the smoothness of the evaluation map 
of $C^\infty(M,\k)$, verified in (1). 

Now we assume that $K$ is regular and put $\g := C^\infty(M,\k)$ and 
$G := C^\infty(M,K)$. 
Then we obtain for each $\xi \in C^\infty(I,\g) \cong C^\infty(I \times M,\k)$ 
(Lemma~A.3) a curve $\gamma \: I \to G$ by 
$\gamma(t)(m) := \Evol_K(\xi^m)(t)$, 
defining a smooth map $I \times M \to K$ (Lemma~A.6(3)), hence a smooth curve in 
$G$ (Lemma~A.2). Now $\delta(\gamma^m) = \xi^m$ implies that 
the evolution map of $G$ is given by 
$\evol_G(\xi)(m) := \evol_K(\xi^m)$.  
Therefore the smoothness of $\evol_G$ follows from 
Lemma~A.3 and the smoothness of the map 
$(\xi,m) \mapsto \evol_K(\xi^m)$ (Lemma~A.6(3)). 

(2) For the Lie group structure we refer to [Wo06]. The holomorphy 
of the evaluation map follows as in (1) from Proposition~I.2(2). 
\qed

\subheadline{Smooth maps with values in regular Lie groups} 

\Definition I.4. Let $M$ be a smooth manifold (with boundary) 
and $K$ a Lie group with Lie algebra~$\k$ 
and Maurer--Cartan form $\kappa_K \in \Omega^1(K,\k)$. 
For an element $f \in C^\infty(M,K)$ we call 
$\delta(f) := f^*\kappa_K =: f^{-1}.df\in \Omega^1(M,\k)$ 
the {\it (left) logarithmic derivative of $f$}. This is a 
$\k$-valued $1$-form on $M$. We thus obtain a map 
$$ \delta \: C^\infty(M,K) \to \Omega^1(M,\k) $$
satisfying the cocycle condition 
$$ \delta(f_1 f_2) = \Ad(f_2)^{-1}.\delta(f_1) + \delta(f_2).\leqno(1.2) $$
(cf.~[KM97, 38.1], [GN07]). From this it easily follows that if $M$ is connected, then 
$$\delta(f_1) = \delta(f_2) \quad \Longleftrightarrow \quad 
(\exists k \in K)\ f_2 = \lambda_k \circ f_1. \leqno(1.3) $$
If $K$ is abelian, then $\delta$ 
is a group homomorphism whose kernel consists of the locally constant maps~$M \to K$. 

We call $\alpha \in \Omega^1(M,\k)$ {\it integrable} if there exists a smooth function 
$f \: M \to K$ with $\delta(f) = \alpha$. We say that $\alpha$ is {\it locally 
integrable} if each point $m \in M$ has an open neighborhood $U$ such that 
$\alpha\res_U$ is integrable. We note that for any smooth map 
$f \: M\to K$ and any smooth curve $\gamma \: [0,1] \to M$ we have 
$$ f(\gamma(1)) 
= f(\gamma(0)) \evol_K(\gamma^*\delta(f)) 
= f(\gamma(0)) \evol_K(\delta(f \circ \gamma)). \leqno(1.4) $$ 
\qed

To describe necessary conditions for the integrability of an element 
$\alpha \in \Omega^1(M,\k)$, we 
define for a manifold $M$ and a locally convex Lie
algebra $\k$, the bracket 
$$ [\cdot, \cdot] \: \Omega^1(M,\k) \times \Omega^1(M,\k) \to \Omega^2(M,\k) $$
by 
$$ \eqalign{ [\alpha, \beta]_p(v,w) 
&:= [\alpha_p(v), \beta_p(w)] - [\alpha_p(w), \beta_p(v)] 
\quad \hbox{ for } \quad v,w \in T_p(M). \cr} $$
Note that $[\alpha, \beta] = [\beta, \alpha]$. 
For a locally convex Lie algebra $\k$ and a smooth manifold $M$ (with boundary), we write 
$$ \MC(M,\k) := \big\{ \alpha \in \Omega^1(M,\k) \: d\alpha 
+ {\textstyle{1\over 2}} [\alpha,\alpha] = 0\big\} $$
for the set of solutions of the {\it Maurer--Cartan equation}. 

The following theorem characterizes the image of 
$\delta$ for a regular Lie group~$K$. 

\Theorem I.5. {\rm(Fundamental Theorem for Lie group-valued functions)} Let 
$K$ be a regular Lie group and $\alpha \in \Omega^1(M,\k)$. 
\litem{(1)} $\alpha$ is locally integrable if and only if $\alpha \in \MC(M,\k)$. 
\litem{(2)} If $M$ is $1$-connected and $\alpha$ is locally integrable, then it is 
integrable. 
\litem{(3)} Suppose that $M$ is connected, fix $m_0 \in M$ and let 
$\alpha \in \MC(M,\k)$. Using piecewise smooth representatives 
of homotopy classes, 
we obtain a well-defined group homomorphism 
$$ \per^{m_0}_\alpha \: 
\pi_1(M,m_0) \to K, \quad [\gamma] \mapsto \evol_K(\gamma^*\alpha), $$
and $\alpha$ is integrable if and only if this homomorphism is trivial. 

\Proof. (1) and (2) follow directly from [KM97, Th.~40.2] (see also [GN07]). 

(3) (cf.\ [GN07]) 
Let $q_M \: \tilde M \to M$ denote a simply connected covering manifold of $M$ 
 and choose a base point $\tilde m_0 \in \tilde M$ with 
$q_M(\tilde m_0) = m_0$. Then the $\k$-valued $1$-form 
$q_M^*\alpha$ on $\tilde M$ also satisfies the Maurer--Cartan 
equation, so that (2) implies the existence of a unique 
smooth function $\tilde f \: \tilde M \to K$ with 
$\delta(\tilde f) = q_M^*\alpha$ and $\tilde f(\tilde m_0) = \1$. 

We write 
$$ \sigma \: \pi_1(M, m_0) \times 
\tilde M \to \tilde M, \quad (d,m) \mapsto d.m  =: \sigma_d(m) $$
for the left 
action of the fundamental group $\pi_1(M,m_0)$ on $\tilde M$. 
In view of (1.3) in Definition~I.4, the relation 
$\delta(\tilde f \circ \sigma_d) = \sigma_d^*q_M^*\alpha = q_M^*\alpha
= \delta(\tilde f)$ 
for each $d \in \pi_1(M,m_0)$ implies the existence of a function 
$$ \chi\: \pi_1(M,m_0) \to K
\quad \hbox{ with } \quad 
\tilde f \circ \sigma_d = \chi(d) \cdot \tilde f\quad \hbox{ for }\ 
d \in \pi_1(M,m_0). $$
For $d_1, d_2 \in \pi_1(M,m_0)$ we then have 
$$ \tilde f \circ \sigma_{d_1 d_2} 
= \tilde f \circ \sigma_{d_1} \circ \sigma_{d_2} 
= (\chi(d_1) \cdot \tilde f) \circ \sigma_{d_2} 
= \chi(d_1) \cdot (\tilde f \circ \sigma_{d_2}) 
= \chi(d_1)\chi(d_2) \cdot \tilde f,$$
so that $\chi$ is a group homomorphism. 
We now pick a smooth lift $\tilde \gamma \: I \to \tilde M$ 
with $q_M \circ \tilde\gamma = \gamma$ and observe that 
$$ \delta(\tilde f \circ \tilde\gamma) 
= \tilde\gamma^*q_M^*\alpha = \gamma^*\alpha, $$
which leads to 
$\chi([\gamma]) = \tilde f([\gamma].\tilde m_0) 
=  \tilde f(\tilde\gamma(1)) = \evol_K(\gamma^*\alpha).$
This proves that $\per^{m_0}_\alpha$ is well-defined and a group homomorphism. 

Clearly, $\per^{m_0}_\alpha$ vanishes if and only if the function $\tilde f$ 
is invariant under the action of $\pi_1(M,m_0)$, which is equivalent to the 
existence of a smooth function $f \: M \to K$ with $f \circ q_M = \tilde f$, 
which in turn means that $\alpha$ is integrable. 
\qed

\Remark I.6. Let $M$ be a connected smooth manifold (with boundary), $m_0 \in M$, and 
$q_M \: \tilde M \to M$ a universal covering map. Further let $K$ be a regular Lie 
group with Lie algebra~$\k$. Then we have an embedding 
$$ q_M^* \: \MC(M,\k) \to \MC(\tilde M,\k)^{\pi_1(M,m_0)}, $$
where the right hand side denotes the set of all solutions of the 
Maurer--Cartan equation on $\tilde M$ which are invariant under the action 
of the fundamental group $\pi_1(M,m_0)$ by deck transformations. 

(a) If $f \in C^\infty(\tilde M, K)$ satisfies $\delta(f) = q_M^*\alpha$, then 
$$ f(d.x) = \per^{m_0}_\alpha(d) \cdot f(x) \quad \hbox{ for all } \quad x \in \tilde M. $$
If, conversely, $f \in C^\infty(\tilde M, K)$ satisfies 
$f(d.x) = \chi(d)\cdot f(x)$ for $x \in \tilde M$ and $d \in \pi_1(M,m_0)$, 
then $\delta(f)$ is $\pi_1(M,m_0)$-invariant, hence of the form 
$\delta(f) = q_M^*\alpha$, and $\chi = \per^{m_0}_\alpha$. This shows that 
$$ C^\infty(\tilde M,K)^\sharp := \{ f \in C^\infty(\tilde M,K) \: 
\delta(f) \in q_M^* \MC(M,\k)\} $$
is fibered over the set $\Hom(\pi_1(M,m_0), K)$ by 
$$ C^\infty(\tilde M,K)^\sharp
= \bigcup_{\chi \in \Hom(\pi_1(M,m_0), K)} C^\infty(\tilde M,K)_\chi, $$
where 
$$ C^\infty(\tilde M,K)_\chi 
:= \{ f \in C^\infty(\tilde M,K)\: (\forall x \in \tilde M)(\forall d \in \pi_1(M,m_0))\ 
f(d.x) = \chi(d)f(x)\}. $$
The set $C^\infty(\tilde M,K)_\chi$ is invariant under multiplication 
with functions in $q_M^*C^\infty(M,K)$ from the right. Conversely, 
for $f, g \in C^\infty(\tilde M,K)_\chi$, the function 
$g^{-1} \cdot f$ factors through a function $F$ on $M$ with 
$g \cdot q_M^*F = f$. Therefore the fibers of 
$C^\infty(\tilde M,K)^\sharp$ coincide with the orbits of the 
group $C^\infty(M,K)$, acting by right multiplication. 
On the set $\MC(M,\k)$, the corresponding action is given by 
$$ \alpha * f := \delta(f)+ \Ad(f)^{-1}.\alpha $$
for $\alpha \in \MC(M,\k)$ and $f \in C^\infty(M,\k)$ (Definition~I.4). 

In general, not all the sets $C^\infty(\tilde M,K)_\chi$ are non-empty. 
Indeed, this condition can be interpreted as the smooth (which is 
equivalent to the topological [MW06]) 
triviality of the flat principal $K$-bundle 
$P_{\chi}:=\widetilde{M}\times_{\chi}K := (\tilde M \times K)/\pi_1(M)$, where 
$\pi_1(M)$ acts on $\tilde M \times K$ by $d.(x,k) = (d.x, \chi(d)k)$. 
Not all such bundles are 
topologically  
trivial. If $\Sigma$ is a compact Riemann surface of genus $g > 1$, then there 
are flat non-trivial $\SL_2(\R)$-bundles over $\Sigma$ 
(see p.~24 in [KT68], which uses results from [Mil58]). 
On the other hand, if $K$ is a complex algebraic group, then 
[Gro68, Cor.~7.2] asserts that all rational characteristic classes of the 
bundles $P_\chi$ vanish.  
 
(b) Note that for 
$f_i \in C^\infty(\tilde M,K)_{\chi_i}$, $i =1,2$, and 
$\im(\chi_2) \subeq Z(K),$
we have $\chi_1\chi_2 \in \Hom(\pi_1(M,m_0), K)$ with 
$f_1 f_2 \in C^\infty(\tilde M,K)_{\chi_1\chi_2}.$
For $\delta(f_i) = q_M^*\alpha_i$ we then have $\delta(f_1f_2) = q_M^*\beta$, where 
$$ \beta = \alpha_2 + \oline{\Ad(f_2)}^{-1}.\alpha_1 $$
and $\oline{\Ad(f_2)}$ is the well-defined function $M \to \Aut(\k)$, defined by 
$\oline{\Ad(f_2)} \circ q_M = \Ad(f_2)$ (cf.\ Definition~I.4). 
\qed

For later use in Section III, we record the following formula: 

\Lemma I.7. Let $x \in \k$ and $\beta \in \Omega^1(M,\R)$ be a closed $1$-form. 
Then $\alpha := \beta \cdot x \in \Omega^1(M,\k)$ satisfies 
$$ \per^{m_0}_\alpha = \exp_K \circ (\per^{m_0}_\beta \cdot x). \leqno(1.5) $$

\Proof.  That $\alpha$ satisfies the Maurer--Cartan equation follows from 
$d\alpha = d\beta \cdot x = 0 = [\alpha,\alpha]$. To calculate $\per^{m_0}_\alpha$, 
we first pick a smooth function $h \in C^\infty(\tilde M,\R)$ with 
$dh = q_M^*\beta$. Then 
$$ f \: \tilde M \to K, \quad m \mapsto \exp_K(h(m)x) $$
is a smooth function with 
$$ \eqalign{ T_m(f)v 
&= T_{h(m)x}(\exp_K)(dh(m)v \cdot x)
= dh(m)v \cdot T_{h(m)x}(\exp_K)(x)\cr
&= dh(m)v \cdot T_\1(\lambda_{\exp_K(h(m)x)}) x 
= dh(m)v \cdot T_\1(\lambda_{f(m)}) x, \cr} $$
so that 
$$ \delta(f) = dh \cdot x = q_M^*(\beta \cdot x) = q_M^*\alpha. $$
Let $\tilde m_0 \in q_M^{-1}(m_0)$ and assume w.l.o.g.\ that 
$h(\tilde m_0) = 0$, so that $f(\tilde m_0) = \1$. We then have 
$$ \per^{m_0}_\alpha(d) = f(d.\tilde m_0) = \exp_K(h(d.\tilde m_0)x) 
= \exp_K(\per^{m_0}_\beta(d) x), $$
which implies (1.5). 
\qed

\subheadline{Uniqueness of regular Lie group structures} 

The following theorem characterizes smooth maps $N \times M \to K$ in terms 
of smoothness of maps defined on $N$. We shall need it later to prove 
the uniqueness of a regular Lie group structure on $C^\infty(M,K)$ with 
smooth evaluation map. 

\Proposition I.8. Let $N$ be a locally convex manifold, 
$M$ a connected finite-dimensional manifold and $K$ a regular Lie group. 
Then a function $f\: N \times M \to K$ is smooth if and only if 
\litem{(1)} there exists a point $m_0 \in M$ for which the map 
$f^{m_0} \: N \to K, n \mapsto f(n,m_0)$ is smooth, and 
\litem{(2)} the functions $f_n \: M \to K, m \mapsto f(n,m)$ are smooth and 
$F \: N \to \Omega^1(M,\k), n \mapsto \delta(f_n)$ is smooth. 

\Proof. ``$\Rarrow$'': If $f$ is smooth, then all the maps $f^m$ and $f_n$ are smooth. 
To see that $F$ is smooth, we recall that $\Omega^1(M,\k)$ is a closed subspace of 
$C^\infty(TM,\k)$ (Definition~I.1(d)), 
so that it suffices to show that the map 
$$ \tilde F \: N \times  TM \to \k, \quad 
(n,v) \mapsto \delta(f_n)v = \kappa_K(T(f_n)v) = \kappa_K(T(f)(0,v)) $$
is smooth (Lemma~A.2). Since the Maurer--Cartan form $\kappa_K$ of $K$ is a
smooth map $TK \to \k$, the assertion follows from the smoothness of 
$T(f) \: T(N \times M) \cong TN \times TM \to TK$. 

``$\Larrow$'': {\bf Step 1:} First we show that $f^m$ is smooth for each 
$m \in M$. Pick a smooth path $\gamma \: [0,1] \to M$ with 
$\gamma(0)= m_0$ and $\gamma(1) = m$. Then 
$$ f^m(n) = f_n(m) 
= f_n(m_0) \evol_K(\delta(f_n \circ \gamma)) 
= f_n(m_0) \evol_K(\gamma^*\delta(f_n))
= f^{m_0}(n) \evol_K(\gamma^*F(n)). $$
Hence the smoothness of $f^{m_0}$, $\evol_K$, 
$F$,  and the continuity of the linear map 
$\gamma^* \: \Omega^1(M,\k) \to C^\infty([0,1],\k)$ imply that $f^m$ is smooth. 

{\bf Step 2:} Now we show that $f$ is smooth. To this end, let 
$m \in M$ and choose a chart $(\phi, U)$ of $M$ for which 
$\phi(U)$ is convex with $\phi(m) = 0$. We have to show that the map 
$$ h \: N \times U \to K, \quad (n,x) \mapsto f(n,\phi^{-1}(x)) $$
is smooth. For $\gamma_x(t) := tx$, $0 \leq t \leq 1$, we have 
$$ h(n,x) = h(n,\gamma_x(1)) 
= h(n,0) \evol_K(\delta(f_n \circ \phi^{-1} \circ \gamma_x))
= f^m(n) \evol_K(\gamma_x^* (\phi^{-1})^* F(n)). $$
Since $f^m$ and $\evol_K$ are smooth and 
$$ (\phi^{-1})^* \: \Omega^1(U,\k) \to \Omega^1(\phi(U),\k)  $$
is a topological linear isomorphism, it suffices to show that the map 
$$ \Omega^1(\phi(U),\k) \times U \to C^\infty([0,1],\k), \quad 
(\alpha,x) \mapsto \gamma_x^*\alpha $$
is smooth. In view of Lemma~A.2, this follows from the smoothness 
of the map 
$$ \Omega^1(\phi(U),\k) \times U \times [0,1] \to \k, \quad 
(\alpha,x,t) \mapsto (\gamma_x^*\alpha)_t = \alpha_{tx}(x), $$
which is a consequence of the smoothness of the evaluation map of 
$C^\infty(TM,\k)$ (Proposition~I.2). 
\qed

The following theorem characterizes the ``good'' Lie group structures on 
$C^\infty(M,K)$ in various ways. Note that we do not assume that the Lie group 
structure is compatible with the topology on $C^\infty(M,K)$. 

\Proposition I.9. Let 
$M$ be a connected finite-dimensional smooth manifold and $K$ a regular 
Lie group. Suppose that the group 
$G := C^\infty(M,K)$ carries a Lie group structure for 
which $\g := C^\infty(M,\k)$ is the corresponding 
Lie algebra and all evaluation maps 
$\ev_{m} \: G \to K$, $m \in M$, are smooth with 
$$\L(\ev_m) = \ev_m \: \g \to \k. $$
Then the following assertions hold: 
\litem{(1)} The evaluation map $\ev \: G \times M \to K, (f,m) \mapsto f(m)$ is smooth. 
\litem{(2)} If, in addition, $G$ is regular, 
then a map $f \: N \to G$ ($N$ a locally convex smooth manifold) 
is smooth if and only if the corresponding map $f^\wedge \: N \times M \to K$ 
is smooth. 

\Proof. (1) Let $N \subeq M$ be a compact submanifold (possibly with boundary). Then 
$C^\infty(N,K)$ carries the structure of a regular Lie group 
(Proposition~I.3). Let $q_G \: \tilde G_0 \to G_0$ denote the universal covering of 
the identity component $G_0$ of $G$. Consider the continuous homomorphism of 
Lie algebras 
$$ \psi \: \L(G) = C^\infty(M,\k) \to C^\infty(N,\k), \quad f \mapsto f\res_N. $$
In view of the regularity of $C^\infty(N,K)$, there exists a unique morphism 
of Lie groups 
$$ \tilde\phi \: \tilde G_0\to C^\infty(N,K)  
\quad \hbox{ with } \quad \L(\tilde\phi) = \psi, $$
where $\tilde G_0$ is the universal covering group of the identity component $G_0$ 
of $G$.
Then, for each $n \in N$, the homomorphism 
$\ev_n \circ \tilde\phi \: \tilde G_0 \to K$ is smooth 
with differential 
$\L(\ev_n \circ \tilde \phi) = \ev_n,$ 
so that 
$\ev_n \circ \tilde\phi = \ev_n \circ q_G,$
where $q_G \: \tilde G_0 \to G$ 
is the universal covering map. We conclude that 
$$ \ker q_G \subeq \ker \tilde\phi, $$
and hence that $\tilde{\phi}$ factors through the restriction map 
$G_0 = C^\infty(M,K)_0 \to C^\infty(N,K).$
In particular, the restriction map 
$C^\infty(M,K) \to C^\infty(N,K)$ is a smooth homomorphism of Lie groups. 

This implies in particular that for each relatively compact open subset 
$U \subeq M$, the map 
$G \to \Omega^1(U,\k), f \mapsto \delta(f\res_U)$ is smooth (cf.\ Lemma~A.5(1)), 
and since $\Omega^1(M,\k)$ embeds into $\prod_U \Omega^1(U,\k)$, it follows that 
$\delta$ is smooth. 

(2) Now we assume that the evaluation map is smooth. If $f$ is smooth, then 
$f^{m_0} = \ev_{m_0} \circ f$ is smooth and (1) entails that 
$\delta \circ f \: N \to \Omega^1(M,\k)$ is smooth, so that 
Proposition~I.8 implies that $f^\wedge$ is smooth. 

If, conversely, $f^\wedge$ is smooth, we have to show that 
$f$ is smooth. To this end, we may w.l.o.g.\ assume that $N$ 
is $1$-connected, since the assertion is local with respect to 
$N$. We define $\beta \in \Omega^1(N,\g)$ by 
$$ \beta v = \kappa_K(T(f^\wedge)(v,0)), $$
which shows immediately that $\beta$ defines a smooth 
map $TN \times M \to \k$ which is linear on the tangent spaces of $N$, 
and with 
$C^\infty(TN \times M,\k) \cong C^\infty(TN,\g)$ (Lemma A.3), we see that 
this is an element of $\Omega^1(N,\g)$. 

We claim that $\beta$ satisfies the Maurer--Cartan equation. 
In fact, for each $m \in M$, we have 
$$ \ev_m \circ \beta = \delta(\ev_m \circ f^\wedge), $$
which satisfies the  Maurer--Cartan equation. Since 
the evaluation map $\ev_m \: \g \to \k$ is a homomorphism of 
Lie algebras, and the corresponding maps 
$\ev_m \: \Omega^1(M,\g) \to \Omega^1(M,\k)$ separate the points, 
it follows that $\beta$ satisfies the Maurer--Cartan equation. 

Fix a point $n_0 \in N$. Since $G$ is regular and $N$ is $1$-connected, 
the Fundamental Theorem (Theorem~I.5) implies the existence of a 
unique smooth function $h \: N \to G$ with 
 $h(n_0) = f(n_0)$ and $\delta(h) = \beta$. For each 
$m \in M$ we then have 
$h(n_0)(m) = f(n_0)(m)$ and 
$$ \delta(\ev_m \circ h) 
= \ev_m \circ \delta(h) 
= \ev_m \circ \beta = \delta(\ev_m \circ f), $$
so that the uniqueness part of the Fundamental Theorem, applied to $K$-valued 
functions, yields $\ev_m \circ f = \ev_m \circ h$ for each $m$, which leads to 
$h = f$. This proves that $f$ is smooth. 
\qed

In the following, we say that a Lie group structure on $C^\infty(M,K)$ 
is {\it compatible with evaluations} if it satisfies the 
assumptions of the preceding proposition. 
The following corollary contains one of the main results of this section. 
It asserts that there is at most one regular Lie group structure 
compatible with evaluations. 

\Corollary I.10. Under the assumptions of the preceding theorem, 
there exists at most one regular Lie group structure 
 on the group $G := C^\infty(M,K)$ compatible with evaluations. 

\Proof. Let $G_1$ and $G_2$ be two regular Lie groups 
obtained from Lie group structures on $C^\infty(M,K)$ compatible 
with evaluations. In view of Proposition~I.9, 
the smoothness of the evaluation maps 
$\ev_j \: G_j \times M \to K$ implies that the identity maps 
$G_1 \to G_2$ and $G_2 \to G_1$ are smooth, hence that 
$G_1$ and $G_2$ are isomorphic Lie groups. 
\qed

In the preceding corollary, we do not have to assume that the Lie group 
structure on $C^\infty(M,K)$ is compatible with the smooth compact open 
topology, but even if we consider only regular Lie group structures 
compatible with this topology, the uniqueness of these structures 
does not directly follow from Lie theoretic considerations: 

\Remark I.11. If $G$ is a regular Lie group and $H$ any simply connected Lie group, 
then any continuous homomorphism $\psi \: \L(H) \to \L(G)$ of Lie algebras 
integrates to a unique homomorphism $\phi \: H \to G$ with 
$\L(\phi) = \psi$ (cf.\ [Mil84]). This implies in particular that 
two $1$-connected regular Lie groups with isomorphic Lie algebras are isomorphic. 
In this sense regularity of a Lie group is crucial for uniqueness results. 

On the other hand, we do not know if there exist topologically isomorphic 
regular Lie groups $G_1$ and $G_2$ which are not isomorphic as Lie groups. 
To prove that this is not the case, 
we would need a result on the automatic smoothness of continuous 
homomorphisms of Lie groups, but presently the optimal result in this direction 
requires at least H\"older continuity (cf.\ [Gl05]). 
\qed

\sectionheadline{II. Logarithmic derivatives and the Maurer--Cartan equation} 

\nin In this section we describe a strategy to 
obtain a Lie group structure on the group 
$G := C^\infty(M,K)$ for a non-compact connected manifold $M$. It is based  
on the injectivity of the logarithmic derivative on the normal subgroup 
$$G_* := \{ f \in C^\infty(M,K) \: f(m_0) = \1\},$$ 
where $m_0$ is a base point in $M$. We thus realize 
$G_*$ as a subset of $\MC(M,\k) \subeq \Omega^1(M,\k)$. 
There are several situations, in which one can show that $\delta(G_*)$ 
is a manifold, and where transferring the manifold structure of $\im(\delta)$ to 
$G_*$ leads to a Lie group structure on $G_*$ and hence on 
$G \cong G_* \rtimes K$ because $K$ (realized as the constant functions on $M$)  
acts smoothly by conjugation on 
$\Omega^1(M,\k) \supeq \delta(G_*)$ (Lemma~A.5(2)). 

One of the main results of this section 
is Theorem~II.2, asserting that $\delta(G_*)$ is a Lie group whenever it is a 
submanifold of $\Omega^1(M,\k)$. This condition is satisfied in particular 
if $M$ is one-dimensional (Corollary~II.3). Then we establish an iterative 
procedure leading to regular Lie group structures on 
$C^\infty(\R^n \times M,K)$ for any compact smooth manifold $M$ and any~$n$ 
(Corollary~II.5).

\Proposition II.1. If $K$ is a regular Lie group and $M$ is 
a connected finite-dimensional smooth 
manifold, then the map 
$$ \delta \: C^\infty_*(M,K) \to \Omega^1(M,\k) $$
is a topological embedding. Let $\Evol_K := \delta^{-1} \: \im(\delta) \to C^\infty_*(M,K)$ 
denote its inverse. Then $\delta$ is an isomorphism of 
topological groups if we endow $\im(\delta)$ with the group structure defined by 
$$ \alpha * \beta := \beta + \Ad(\Evol_K(\beta))^{-1}.\alpha \leqno(2.1) $$
and 
$$ \alpha^{-1} = - \Ad(\Evol_K(\alpha)).\alpha, \leqno(2.2)$$

\Proof. First we show that $\delta$ is continuous. By definition of the topology 
on $C^\infty(M,K)$, the tangent map induces a continuous group homomorphism 
$$ T \: C^\infty(M,K) \to C^\infty(TM,TK), \quad f \mapsto T(f). $$
Let $\kappa_K \: TK \to \k$ denote the (left) Maurer--Cartan form of $K$.
Since $\delta(f) = f^*\kappa_K = \kappa_K \circ T(f)$, it follows that the composition 
$$ C^\infty(M,K) \to C^\infty(T(M), T(K)) \to C^\infty(T(M),\k), \quad 
f \mapsto T(f) \mapsto \delta(f) $$
is continuous. 

Next we show that $\delta$ is an embedding. Consider 
$\alpha = \delta(f)$ with $f \in C^\infty_*(M,K)$, i.e., 
$f(m_0) = \1$ holds for the base point $m_0 \in M$. 
To reconstruct $f$ from $\alpha$, we pick for $m \in M$ a piecewise 
smooth path $\gamma \: [0,1] \to M$ with $\gamma(0) = m_0$ and $\gamma(1) = m$. 
Then $\delta(f\circ \gamma) = \gamma^*\delta(f) = \gamma^*\alpha$ implies 
$f(m) = \evol_K(\gamma^*\alpha)$. 

We now choose an open neighborhood 
$U$ of $m$ and a chart $(\phi,U)$ of $M$ such that 
$\phi(U)$ is convex with $\phi(m) = 0$. Then, for each $x \in U$, 
(1.4) in Definition~I.4 yields 
$$ f(x) = f(m) \cdot \evol_K(\gamma_x^*\alpha), \leqno(2.3) $$
where $\gamma_x(t) = \phi^{-1}(t\phi(x))$. 

>From Lemma~A.6(1),(2), we immediately derive that the map 
$$ \Omega^1(M,\k) \times U \to K, \quad 
(\alpha,x) \mapsto \evol_K(\gamma^*\alpha) \cdot \evol_K(\gamma_x^*\alpha)\leqno(2.4)  $$
is smooth, so that the corresponding map 
$\Omega^1(M,\k) \to C^\infty(U,K)$
is in particular continuous 
(Lemma~A.1). We conclude that the map 
$$ \delta(C^\infty_*(M,K)) \to C^\infty(U,K), \quad \delta(f) \mapsto f\res_U $$
is continuous. We finally observe that for each open covering 
$M = \bigcup_{j \in J} U_j$, 
the restriction maps to $U_i$ lead to a topological 
embedding $C^\infty(M,K) \into \prod_{j \in J} C^\infty(U_j,K)$, 
and this 
completes the proof. 
\qed

\Theorem II.2. Let $M$ be a connected finite-dimensional smooth 
manifold (with boundary) and $K$ a regular Lie group. 
Assume that the subset $\delta(C^\infty_*(M,K))$ is a smooth submanifold of 
$\Omega^1(M,\k)$ and endow 
$C^\infty_*(M,K)$ with the manifold structure for which 
$\delta \: C^\infty_*(M,K) \to \im(\delta)$
is a diffeomorphism and 
$$C^\infty(M,K) \cong K \ltimes C^\infty_*(M,K)$$
with the product manifold structure. 
Then the following assertions hold: 
\litem{(1)} For each locally convex manifold $N$, 
a map $f \: N \times M \to K$ is smooth if and only 
if all the maps $f_n \: M \to K, m \mapsto f(n,m)$ are smooth 
and the corresponding map 
$$ f^{\vee} \: N \to C^\infty(M,K), \quad n \mapsto f_n $$
is smooth. 
\litem{(2)} $K$ acts smoothly by conjugation on $C^\infty_*(M,K)$, and 
$C^\infty(M,K)$ carries a regular Lie group structure compatible with evaluations. 

\Proof. (1) Let $m_0$ be the base point of $M$. 
According to Proposition~I.8, $f \: N \times M \to K$ is smooth if and only 
if $f^{m_0}$ is smooth, all the maps $f_n$ are smooth, 
and $\delta \circ f^\vee \: N \to \Omega^1(M,\k)$ is smooth. 
In view of our definition of the manifold structure on 
$C^\infty_*(M,K)$, the latter condition is equivalent to the smoothness 
of the map $N \to C^\infty_*(M,K), n \mapsto f_n(m_0)^{-1}f_n = f^{m_0}(n)^{-1} f_n$. 
Since the evaluation in $m_0$ coincides with the projection 
$$ C^\infty(M,K) \cong K \ltimes C^\infty_*(M,K) \to K, $$
we see that $f$ is smooth if and only if 
all the maps $f_n$ are smooth and $f^\vee$ is smooth. 

(2) For the evaluation map $f = \ev \: G \times M \to K$, 
we have $\ev^\vee = \id_G$ and $\ev_g = g$ for each $g \in G$. Hence 
(1) implies that $\ev$ is smooth. 

In view of Proposition~II.1, $\delta$ is an isomorphism of topological 
groups if $\im(\delta)$ is endowed with the group structure (2.1). 
We now show that the operations (2.1) and (2.2) are smooth 
with respect to the submanifold structure on $\im(\delta)$. 

{\bf The Lie group structure:} It suffices to show that the map 
$$ \im(\delta) \times \Omega^1(M,\k) \to \Omega^1(M,\k), \quad 
(\alpha,\beta) \mapsto \Ad(\Evol_K(\alpha)).\beta $$
is smooth. For each open covering $(U_j)_{j \in J}$, we obtain an embedding 
$\Omega^1(M,\k) \into \prod_{j \in J} \Omega^1(U_j,\k)$, so that it 
suffices to prove for each $m \in M$ the existence of an open 
neighborhood $U$ of $m$, for which the map 
$$ \im(\delta) \times \Omega^1(U,\k) \to \Omega^1(U,\k), \quad 
(\alpha,\beta) \mapsto \Ad(\Evol_K(\alpha)).\beta $$
is smooth. Choosing $U$ so small that it lies in a chart domain, 
we have $\Omega^1(U,\k) \cong C^\infty(U,\k)^d$ for $d = \dim M$, so that it 
suffices to show that 
$$ \im(\delta) \times C^\infty(U,\k) \to C^\infty(U,\k), \quad 
(\alpha,f) \mapsto \Ad(\Evol_K(\alpha)).f, $$
is smooth. Now it suffices to see that the map 
$$ \im(\delta)  \times C^\infty(U,\k) \times U \to \k, \quad 
(\alpha,f,x) \mapsto \Ad(\Evol_K(\alpha)(x)).f(x) $$
is smooth. 
Since the action map $K \times C^\infty(M,\k) \to C^\infty(M,\k)$ 
is smooth (Lemma~A.5(1)), it suffices to recall from Proposition~I.2 that the 
evaluation map of $C^\infty(U,\k)$ is smooth and to show that the map 
$$ \im(\delta)  \times U \to K, \quad 
(\alpha,x) \mapsto \Evol_K(\alpha)(x) \leqno(2.5) $$
is smooth. 

Let $m \in M$. To obtain $\Evol_K(\alpha)$, we pick a piecewise 
smooth path $\gamma \: [0,1] \to M$ with $\gamma(0) = m_0$ and $\gamma(1) = m$. 
Then 
$\delta(\Evol_K(\alpha)\circ \gamma) = \gamma^*\delta(\Evol_K(\alpha)) = \gamma^*\alpha$ 
implies 
$$ \Evol_K(\alpha)(m) = \evol_K(\gamma^*\alpha). $$
We now choose an open neighborhood 
$U$ of $m$ and a chart $(\phi,U)$ of $M$ such that 
$\phi(U)$ is convex. Then, for each $x \in U$, (1.4) in Definition~I.4  
entails 
$$ \Evol_K(\alpha)(x) = \Evol_K(\alpha)(m) \cdot \evol_K(\gamma_x^*\alpha),  $$
where $\gamma_x(t) = \phi^{-1}(t\phi(x))$. 
>From Lemma~A.6(1),(2), we derive that the map 
$$ \Omega^1(M,\k) \times U \to K, \quad 
(\alpha,x) \mapsto \evol_K(\gamma^*\alpha) \cdot \evol_K(\gamma_x^*\alpha)  $$
is smooth, so that restriction to the submanifold 
$\im(\delta)$ implies the smoothness of (2.5). 
We conclude that multiplication and the inversion in $\im(\delta)$ 
is smooth and hence that it is a Lie group. 

{\bf The Lie algebra:} Next we verify regularity. 
To this end, we first determine the tangent space $T_0(\im(\delta))$ to 
see the Lie algebra of this group. 
Let $\eta \: I \to \im(\delta)$ be a smooth curve 
with $\eta(0) =0$ and $\beta := \eta'(0)$. Then 
$$ d \eta(t) + {1\over 2}[\eta(t), \eta(t)] = 0 $$
for each $t \in I$ yields 
$ d \eta'(0) = 0$, so that $\beta$ is closed. We also have 
$$ \1 = \per_{\eta(t)}^{m_0}(\gamma) = \evol_K(\gamma^*\eta(t)) $$ 
for each smooth loop $\gamma$ in $m_0$ and each $t \in I$. 
Taking the derivative in $t = 0$, we get with Lemma~A.5(1): 
$$ 0 = T_0(\evol_K)(\gamma^*\beta) = \int_0^1 \gamma^*\beta = \int_\gamma \beta. $$ 
Hence all periods of $\beta$ vanish, so that $\beta$ is exact. If, 
conversely, $\beta \in \Omega^1(M,\k)$ is an exact $1$-form, then 
$\beta = df$ for some $f \in C^\infty_*(M,\k)$, and the curve 
$\alpha(t) := \delta(\exp_K(tf))$ in $\im(\delta)$ satisfies 
$\alpha'(0) = T_\1(\delta)f = df = \beta$. This shows that 
$$ T_0(\im(\delta)) = B^1_{\rm dR}(M,\k) = dC^\infty_*(M,\k) \cong C^\infty_*(M,\k), $$ 
as a topological vector space (apply Proposition~II.1 to the Lie group $(\k,+)$).  

Next, recall that for each $m \in M$ and $\gamma \: [0,1] \to M$ from 
$m_0$ to $m$ we have 
$$ \Evol_K(\alpha)(m) = \evol_K(\gamma^*\alpha), $$
so that we get for any smooth curve $\eta$ in $\im(\delta)$ 
with $\eta(0) = 0$ and $\eta'(0) = df$ with $f \in C^\infty_*(M,\k)$ the relation 
$$ \derat0 \Evol_K(\eta(t))(m) = \derat0 \evol_K(\gamma^*\eta(t)) 
= T_0(\evol_K)\gamma^*\eta'(0) = \int_0^1 \gamma^*\eta'(0) 
=  f(m), $$
i.e., 
$$ \derat0 \Evol_K(\eta(t)) = f. $$
Now we can determine the Lie bracket on $T_0(\im(\delta))$. 
Let $\eta_j \: I \to \im(\delta)$, $j =1,2$, be smooth curves 
in $\1$ and $f_j \in C^\infty_*(M,\k)$ with $df_j = \eta'(0)$. 
Then 
$$ {\partial^2 \over \partial s \partial t}\mid_{s=t=0} 
\eta_1(s) * \eta_2(t) 
= \derat0 \Ad(\Evol_K(\eta_2(t)))^{-1}.\eta_1'(0) 
= -[f_2, df_1]. $$
For the Lie bracket in $\L(\im(\delta)) = dC^\infty(M,\k)$, we thus obtain the 
formula 
$$ [df_1, df_2] 
= -[f_2, df_1] + [f_1, df_2] = d[f_1, f_2], $$
showing that 
$$ d \: C^\infty_*(M,\k) \to T_0(\im(\delta)) $$
is an isomorphism of Lie algebras if $C^\infty_*(M,\k)$ is endowed with the 
pointwise Lie bracket. 

{\bf Regularity:} 
It remains to verify the regularity of $G$, 
i.e., the smoothness of the map 
$$ \evol_G \: C^\infty(I,\g) \to G \cong C^\infty(M,K).$$ 
First we make this map more explicit. 
Let $\xi \in C^\infty(I,\g) \cong C^\infty(I \times M, \k)$ (Lemma~A.2). 
To see that the curve 
$\gamma(t)(m) := \Evol_K(\xi^m)(t)$ in $G$ is smooth, we observe that 
the map 
$I \times M \to K, (t,m) \mapsto \gamma(t)(m)$ is smooth (Lemma~A.6(3)), 
so that (1) implies that $\gamma \: I \to G$ is smooth. 
We also obtain from Lemma~A.6(3) that $\delta(\gamma)_t = \xi_t$, so that 
$\delta(\gamma)= \xi$. Hence 
$$ \evol_G(\xi)(m) =  \evol_K(\xi^m) = \evol_G^\vee(\xi,m). $$
In view of (1), the smoothness of $\evol_G^\vee$ (Lemma~A.6(3)) 
now implies the smoothness of $\evol_G$. 
\qed

If $M$ is one-dimensional, then each $\k$-valued $2$-form on $M$ is
trivial, so that $d\alpha = 0 = [\alpha,\alpha]$ 
for each $\alpha \in \Omega^1(M,\k)$. Therefore 
all $1$-forms trivially solve the Maurer--Cartan equation. We thus obtain: 

\Corollary II.3. If $M$ is a one-dimensional $1$-connected 
real manifold (with boundary), then the group 
$C^\infty_*(M,K)$ carries a regular Lie group structure 
 for which 
$$ \delta \: C^\infty_*(M,K) \to \Omega^1(M,\k) \cong C^\infty(M,\k)  $$
is a diffeomorphism and $C^\infty(M,K) \cong C^\infty_*(M,K) \rtimes K$ 
carries the structure of a regular Lie group 
compatible with evaluations and the smooth compact open topology. 
\qed

For the case $M = \R$, the preceding corollary can also be found in the 
book of Kriegl and Michor ([KM97, Th.~38.12]). 
Note that any $1$-connected $\sigma$-compact $1$-dimensional manifold 
with boundary is diffeomorphic to 
$\R,$ $[0,1]$ or $[0,\infty[.$

\Lemma II.4. If $K \not= \exp \k$, then 
the exponential image of $C^\infty(\R,\k)$  
is not an identity neighborhood in $C^\infty(\R,K)$. 

\Proof. The exponential function $C^\infty(\R,\k)$ is simply given by 
$\exp(\xi)\,:=\,\exp_K\circ\,\xi$, where $\exp_K$ is the exponential 
function of $K$.  

Let $k \in K \setminus \exp \k$ and consider a smooth curve 
$g \: \R \to K$ with $g(t) = \1$ for $t < 0$ and $g(t) = k$ for $t > 1$. 
Then $g_n(t) := g(t - n)$ defines a sequence in $C^\infty(\R,K)$, 
converging to~$\1$. 
As $g_n(n+1) = k \not\in \exp \k$, none of the curves $g_n$ is contained in 
the image of the exponential function. 
\qed

\subheadline{Iterative constructions} 

Corollary II.3 is much more powerful than it appears at first sight 
because it can be applied inductively to show that for each compact manifold $M$ 
and $k \in \N$ the group $C^\infty(\R^k \times M,K)$ carries a regular 
Lie group structure compatible with evaluations. This result is the main goal 
of this subsection. First we need two lemmas. The more general 
key result is Theorem~II.7. 

\Lemma II.5. Let $(G_n)_{n \in \N}$ be a sequence of Lie groups, 
$\phi_{nm} \: G_m \to G_n$ morphisms of Lie groups 
defining an inverse system, 
$G := \prolim G_n$ the corresponding topological projective limit group and $\phi_n \: G \to G_n$ the canonical maps. 
Assume that $G$ carries a Lie group structure 
with the following properties: 
\litem{(1)} A map $f \: M \to G$ of a smooth manifold $M$ with values in 
$G$ is smooth if and only if all the maps $f_n := \phi_n \circ f$ are smooth. 
\litem{(2)} $\L(G) \cong \prolim \L(G_n)$ as topological Lie algebras, 
with respect to the projective system defined by the morphisms 
$\L(\phi_{nm}) \: \L(G_n) \to \L(G_m)$. 

\nin Then the map 
$$ \Psi \: C^\infty(M,G) \cong \prolim C^\infty(M,G_n), \quad f 
\mapsto (f_n)_{n \in \N} $$
is an isomorphism of topological groups. 

\Proof. First we note that our assumptions imply that 
$$ TG \cong \L(G) \rtimes G 
\cong \prolim \left(\L(G_n) \rtimes G_n \right)
\cong \prolim T(G_n) $$
as topological groups. Moreover, writing $|\L(G)|$ for the topological 
vector space underlying $\L(G)$, considered as an abelian Lie algebra, we have 
$$ \L(TG) \cong |\L(G)| \rtimes \L(G) 
\cong \prolim \left(|\L(G_n)| \rtimes \L(G_n)\right) 
\cong \prolim \L(TG_n), $$
so that the Lie group $TG$ inherits all 
properties assumed for $G$. Hence we may iterate this argument to obtain 
$$ T^k G \cong \prolim T^k G_n $$
for each $k$ and that (1) holds for the Lie group $T^k G$. 

We thus have topological embeddings 
$$ C(T^k M, T^k G)_c  \into \prolim C(T^k M, T^k G_n)_c, $$
which leads to a topological embedding 
$$ C^\infty(M,G) 
\into \prod_{k \in \N_0} C(T^k M, T^k G)_c
\into \prod_{k \in \N_0} \prolim C(T^k M, T^k G_n)_c
\cong \prolim \prod_{k \in \N_0} C(T^k M, T^k G_n)_c, $$
showing that $\Psi$ is a topological isomorphism. 
\qed

\Lemma II.6. If $N$ and $M$ are compact manifolds (possibly with boundary), 
then the map 
$$ \Phi\: C^\infty(N, C^\infty(M,K)) \to C^\infty(N \times M,K), \quad 
f \mapsto f^\wedge $$ 
is an isomorphism of Lie groups. 

\Proof. The bijectivity of $\Phi$ follows from the smoothness of the evaluation 
map of $C^\infty(M,K)$ (Proposition~I.3) and 
Proposition~I.9. To see that $\Phi$ is an isomorphism of Lie groups, 
let $(\phi,U)$ be a $\k$-chart of $K$ with $\phi(\1) = 0$. 
Then $C^\infty(M,U)$ is an open identity neighborhood, so that 
$C^\infty(N, C^\infty(M,U))$ is an open identity neighborhood, 
and so is $C^\infty(N \times M,U)$. That $\Phi$ restricts to a diffeomorphism 
$$ C^\infty(N, C^\infty(M,U))\to C^\infty(N \times M,U) $$ 
now follows from Lemma A.3 which asserts that 
$$ C^\infty(N, C^\infty(M,\k))\to C^\infty(N \times M,\k), \quad 
f \mapsto f^\wedge $$ 
is an isomorphism of topological vector spaces, hence restricts to 
diffeomorphisms on open subsets. 
\qed

\Theorem II.7. Let 
 $K$ be a regular Lie group and $N$ and $M$ finite-dimensional 
smooth $\sigma$-compact manifolds. We assume that 
$G := C^\infty(M,K)$ carries a regular Lie group structure 
compatible with evaluations and the smooth compact open topology. 
If $C^\infty(N,G)$ also carries a 
regular Lie group structure compatible with evaluations
and the smooth compact open topology,  then 
$C^\infty(N \times M,K)$ carries a regular Lie group structure 
compatible with evaluations. Moreover, the canonical map 
$$ \Phi \: C^\infty(N \times M,K) \to C^\infty(N,G), \quad f \mapsto f^\vee $$
is an isomorphism of Lie groups. 

\Proof. In view of Proposition~I.9, the map $\Phi$ is a bijective group 
homomorphism. First we show that it is an isomorphism of topological groups. 

Let $M = \bigcup_n M_n$ be an exhaustion of $M$ 
by compact submanifolds $M_n$ with boundary satisfying $M_n \subeq 
M_{n+1}^0$. Then our definition of the group topology implies that 
$$ G = C^\infty(M,K) \cong \prolim C^\infty(M_n,K)$$ 
as topological groups. Put $G_n := C^\infty(M_n,K)$ and 
recall from Proposition~I.3 that it carries a regular Lie group structure 
compatible with evaluations. 
We also have the isomorphism of topological Lie algebras 
$$\L(G) = C^\infty(M,\k) 
\cong \prolim \L(G_n) \cong \prolim C^\infty(M_n,\k),$$
and Proposition~I.9 implies that we have for each smooth manifold $X$: 
$$ C^\infty(X,G) 
\cong C^\infty(X \times M,K) 
\cong \prolim C^\infty(X \times M_n,K) 
\cong \prolim C^\infty(X, G_n) $$
on the level of groups (without topology). 

Now let $(N_k)_{k \in \N}$ be an exhaustion of $N$ by compact submanifolds 
with boundary. Then Lemmas~II.5 and II.6 lead to the following isomorphisms 
of topological groups: 
$$ \eqalign{ C^\infty(N,G) 
&\cong 
\prolim C^\infty(N,G_n) = \prolim C^\infty(N,C^\infty(M_n,K)) \cr 
&\cong \prolim_k \prolim_n C^\infty(N_k, C^\infty(M_n,K))  
\cong \prolim_k \prolim_n C^\infty(N_k \times M_n,K) \cr 
&\cong C^\infty(N \times M,K). \cr } $$

The preceding isomorphism leads to a regular Lie group structure on the 
topological group $C^\infty(N \times M,K)$. To see that this is the 
unique regular Lie group structure compatible with  evaluations, 
we first observe that all evaluation maps 
$$ \ev_{(n,m)} = \ev_m \circ \ev_n \: C^\infty(N,C^\infty(M,K)) \to K $$
are smooth and then apply Proposition~I.9 to see that the evaluation map 
on $C^\infty(N \times M,K)$ is smooth. 
\qed

Applying Theorem~II.7 and Corollary~II.3 inductively, we obtain: 
\Corollary II.8. Let  $K$ be a regular Lie group, $M$ a finite-dimensional 
compact manifold, $k \in \N_0$ and 
$N := \R^k \times M$. Then 
$C^\infty(N,K)$ carries a regular Lie group structure compatible with evaluations and 
the smooth compact open topology. 
\qed

\Corollary II.9. For each finite-dimensional connected Lie group $M$ and each 
regular Lie group $K$, the group $C^\infty(M,K)$ 
carries a regular Lie group structure for which the evaluation map is smooth. 

\Proof. This follows from the fact that each connected Lie group $M$ is diffeomorphic 
to $\R^n \times C$, where $C$ is a maximal compact subgroup and 
$n = \dim M - \dim C$ ([Ho65]). 
\qed

\sectionheadline{III. The complex case} 

In this section we assume that $M$ is a connected complex manifold without 
boundary of dimension~$d$ and that $K$ is a regular complex Lie group. 

If $E$ is a complex locally convex space, 
we write $\Omega^1_h(M,E)$ for the space of holomorphic $E$-valued $1$-forms on $M$. 
For a complex Lie algebra $\k$, we write 
$$ \MC_h(M,\k) := \MC(M,\k) \cap \Omega^1_h(M,\k) $$
for the set of holomorphic solutions of the Maurer--Cartan equation and topologize 
this space as a subspace of ${\cal O}(TM,\k)$, endowed with the compact open topology. 

In the complex setting, the Fundamental Theorem is easily deduced from the real 
version (Theorem~I.5): 

\Theorem III.1. {\rm(Complex Fundamental Theorem)} Let 
$M$ be a complex manifold and $K$ be a regular complex Lie group. 
\litem{(1)} A smooth function $f \: M \to K$ is holomorphic if and only if 
$\delta(f)$ is a holomorphic $\k$-valued $1$-form. 
\litem{(2)} An element $\alpha \in \Omega^1_h(M,\k)$ is locally integrable 
to a holomorphic  function if and only if it satisfies the Maurer--Cartan equation. 
\litem{(3)} Suppose that $M$ is connected, fix $m_0 \in M$ and assume that 
$\alpha \in \MC_h(M,\k)$. Then 
$\alpha$ is integrable to a holomorphic 
function $M \to K$ if and only if the period homomorphism $\per_\alpha^{m_0}$ is trivial. 

\Proof. (1) If $f$ is holomorphic, then $T(f) \: T(M) \to T(K)$ is holomorphic, 
and since $\kappa_K$ is holomorphic, the same holds for 
$\delta(f) = \kappa_K \circ T(f)$. 

If, conversely, $\delta(f)$ is a holomorphic $1$-form, then each map 
$T_x(f) \: T_x(M) \to T_{f(x)}(K)$ is complex linear, so that 
$f$ is holomorphic. 

(2), (3) In view of (1), this follows from Theorem~I.5. 
\qed

\Remark III.2. Remark I.6(a) carries over to the holomorphic setting as follows. 
For each homomorphism $\chi \: \pi_1(M,m_0) \to K$ we obtain a holomorphic 
flat $K$-bundle $P_\chi := \widetilde{M}\times_{\chi} K$. 
If $M$ is a Stein manifold $M$ and $K$ is Banach, then the Oka principle 
(cf.\ [Rae77, Th.~2.1]) asserts that this 
bundle has a holomorphic section if and only if it has a continuous section. 
This is the case if and only if the corresponding space 
${\cal O}(\tilde M,K)_\chi$ is non-empty. 

A typical example where ${\cal O}(\tilde M,K)_\chi = \eset$ can be obtained 
as follows: Consider the complex manifold $M = \PGL_n(\C) = \PSL_n(\C)$ 
with the universal covering $\tilde M \cong \SL_n(\C)$ and 
identify $\pi_1(M)$ with the cyclic group $C_n := \SL_n(\C) \cap \C^\times \1$ 
of order $n$. Define $\chi \: \pi_1(M) \to \C^\times$ by 
$z = \chi(z)^{-1}\1$. Then 
$$ P_\chi 
= \tilde M \times_\chi \C^\times 
= \SL_n(\C) \times_\chi \C^\times 
\cong \SL_n(\C) \cdot \C^\times \1 = \GL_n(\C), $$
and this $\C^\times$-bundle over $M$ is non-trivial because the corresponding 
surjective homomorphism 
$$ \pi_1(\GL_n(\C)) \cong \Z \to \pi_1(M) \cong C_n $$
does not split since $\Z$ is torsion free.  
\qed

\Proposition III.3. If $M$ is a connected complex manifold without boundary 
and $K$ a regular complex 
Lie group, then the map 
$$ \delta \: {\cal O}_*(M,K) \to \Omega^1_h(M,\k) \subeq {\cal O}(T(M),\k) $$
is a topological embedding if ${\cal O}(M,K)$ carries the compact open 
topology. 

\Proof. To see that the inclusion 
${\cal O}(M,K) \into C^\infty(M,K)$
is a topological embedding for any complex Lie group $K$, it suffices to prove that 
uniform convergence of holomorphic functions in ${\cal O}(M,K)$ implies 
uniform convergence of all tangent maps $T^n(f) \: T^n(M) \to T^n(K)$ on compact 
subsets. 

Suppose that $f_i \to f$ in ${\cal O}(M,K)$ and let 
$C \subeq M$ be a compact subset for which $f(C)$ lies in $U_K$ 
for some holomorphic $\k$-chart $(\phi,U_K)$ of $K$. Then we may w.l.o.g.\ 
assume that $f_i(C) \subeq U_K$ for each $i$. As we have seen in 
Proposition~I.2(2), this implies that, on the interior $C^0$, 
$\phi \circ f_i\res_{C^0} \: C^0 \to \k$ converges to 
$\phi \circ f\res_{C^0}$ in $C^\infty(C^0, \k)$, but this also implies that 
$f_i\res_{C^0} \to f\res_{C^0}$ in $C^\infty(C^0,K)$ (Lemma~A.4). 
Since each compact subset of $M$ can be covered with finitely many 
sets of the form $C^0$, it follows that $f_i \to f$ in $C^\infty(M,K)$. 

Now the assertion of the corollary follows from Proposition~II.1. 
\qed

\Lemma III.4. If $K$ is a regular complex Lie group, then 
$\evol_K \: C^\infty([0,1],\k) \to K$
is holomorphic. 

\Proof. From Corollary~II.3, we know that $C^\infty([0,1],\k)$ 
carries a natural Lie group structure, and for this Lie group structure,  
the map $\evol_K$ is a group homomorphism. In view of the regularity of 
$K$, this map is smooth. To verify its holomorphy, it therefore suffices to 
show that $T_\1(\evol_K)$ is complex linear which is an immediate consequence 
of Lemma~A.5(1). 
\qed

We now adapt Proposition~I.8, Proposition~I.9 and Corollary~I.10 to complex 
Lie groups to derive a complex version of Theorems~II.2. 

\Proposition III.5. Let $N$ be a locally convex complex manifold, 
$M$ a connected fi\-nite-di\-men\-sio\-nal complex manifold and $K$ a regular complex 
Lie group. 
Then a function $f\: N \times M \to K$ is holomorphic if and only if 
\litem{(1)} there exists a point $m_0 \in M$ for which the map 
$f^{m_0} \: N \to K, n \mapsto f(n,m_0)$ is holomorphic, and 
\litem{(2)} the functions $f_n \: M \to K, m \mapsto f(n,m)$ are holomorphic and 
$F \: N \to \Omega^1_h(M,\k), n \mapsto \delta(f_n)$ is holomorphic. 

\Proof. ``$\Rarrow$'' is verified as in the real case (Proposition~I.8). 

``$\Larrow$'': The proof follows the line of the real case. In addition, 
we use that $\evol_K$ is holomorphic (Lemma~III.4) and that 
the pull-back maps $\gamma^*$ are complex linear. 
One point that requires some extra care is the verification of the 
holomorphy of the map 
$$ \Omega^1_h(\phi(U),\k) \times U \to C^\infty(I,\k), \quad 
\gamma_x^*\alpha = \alpha \circ T(\gamma_x). $$
>From the proof of Proposition~I.8 we know that 
it is smooth and it is complex linear in $\alpha$, 
so that its holomorphy follows from its holomorphy in $x$. 
\qed

\Proposition III.6. Let 
$M$ be a connected finite-dimensional complex manifold and $K$ a regular 
complex Lie group. For a complex Lie group structure on the group 
$G := {\cal O}(M,K)$ for which $\g := {\cal O}(M,\k)$ is the corresponding 
Lie algebra and all evaluation maps 
$\ev_{m} \: G \to K$, $m \in M$, are holomorphic with 
$$\L(\ev_m) = \ev_m \: \g \to \k. $$
Then the following assertions hold: 
\litem{(1)} The evaluation map $\ev \: G \times M \to K, (f,m) \mapsto f(m)$ 
is holomorphic. 
\litem{(2)} If, in addition, $G$ is regular, 
then a map $f \: N \to G$ 
is holomorphic if and only if the corresponding map $f^\wedge \: N \times M \to G$ 
is holomorphic. 

\Proof. (1) From Proposition~I.9 it follows that $\delta$ is smooth. 
It also satisfies the cocycle identity 
$$ \delta(f_1 f_2) = \Ad(f_2)^{-1}.\delta(f_1) + \delta(f_2). $$
Since the maps $\Ad(f)$ are complex linear on $\Omega^1_h(M,\k)$, 
it therefore suffices to observe that $T_\1(\delta)(f) = df$ is complex linear in $f$, 
to conclude that $\delta$ is holomorphic. 

(2) If $f$ is holomorphic, then 
$f^{m_0} = \ev_{m_0} \circ f$ is holomorphic and (1) entails that 
$\delta \circ f \: N \to \Omega^1_h(M,\k)$ is holomorphic, so that 
Proposition~III.5 implies that $f^\wedge$ is holomorphic. 

If, conversely, $f^\wedge$ is holomorphic, we first use 
Proposition~I.9 to see that $f$ is smooth. 
That its differential is complex linear 
follows from the holomorphy of $f^\wedge$. 
\qed

In the following, we say that a Lie group structure on ${\cal O}(M,K)$ 
is {\it compatible with evaluations} if it satisfies the 
assumptions of the preceding proposition. 
As in the real case, we obtain: 

\Corollary III.7. Under the assumptions of the preceding theorem, 
there exists at most one regular complex Lie group structure 
 on the group ${\cal O}(M,K)$ which is compatible with evaluations. 
\qed

\Theorem III.8. Let $M$ be a connected complex manifold and $K$ a complex 
regular Lie group. 
Assume that the subset $\delta({\cal O}_*(M,K))$ 
is a complex submanifold of $\Omega^1_h(M,\k)$ 
and endow ${\cal O}_*(M,K)$ with the manifold structure for which 
$\delta \: {\cal O}_*(M,K) \to \im(\delta)$
is biholomorphic and 
$$ {\cal O}(M,K) \cong K \ltimes {\cal O}_*(M,K)$$
with the product manifold structure. 
Then the following assertions hold: 
\litem{(1)} For each locally convex complex manifold $N$, 
a map $f \: N \times M \to K$ is holomorphic if and only 
if all the maps $f_n \: M \to K, m \mapsto f(n,m)$ are holomorphic 
and the corresponding map 
$$ f^{\vee} \: N \to {\cal O}(M,K), \quad n \mapsto f_n $$
is holomorphic. 
\litem{(2)} $K$ acts holomorphically by conjugation on ${\cal O}_*(M,K)$, and 
${\cal O}(M,K)$ carries a regular complex Lie group structure compatible with 
evaluations. 

\Proof. (1) is proved as Theorem~II.2(1). Here we use 
Proposition~III.5 instead of Proposition~I.8. 

(2) For the evaluation map $f = \ev \: G \times M \to K$, 
we have $\ev^\vee = \id_G$ and $\ev_g = g$ for each $g \in G$. Hence 
(1) implies that $\ev$ is holomorphic. 

{\bf The complex Lie group structure:} 
In view of Propositions~II.1 and III.3, 
$\delta$ is an isomorphism of topological 
groups if $\im(\delta)$ is endowed with the group structure (2.1). 
To see that the group operations are holomorphic, we have to show that the map 
$$ \im(\delta) \times \Omega^1_h(M,\k) \to \Omega^1_h(M,\k), \quad 
(\alpha,\beta) \mapsto \Ad(\Evol_K(\alpha)).\beta $$
is holomorphic. Its smoothness has already been verified in 
the proof of Theorem~II.2. Since it is complex linear in $\beta$, it 
is holomorphic in the second argument. 

We claim that it is also holomorphic 
in the first argument $\alpha$. Since the evaluation maps 
$$\ev_m \: \Omega^1_h(M,\k) \to \Hom(T_m(M),\k)$$
are complex linear, and the adjoint action of $K$ on $\k$ is holomorphic, 
it suffices to show that for each element $m \in M$, the map 
$$ \im(\delta) \to K, \quad \alpha \mapsto \Evol_K(\alpha)(m) $$
is holomorphic. 

For any piecewise 
smooth path $\gamma \: [0,1] \to M$ with $\gamma(0) = m_0$ and $\gamma(1) = m$ 
we have 
$$ \Evol_K(\alpha)(m) = \evol_K(\gamma^*\alpha), $$
so that the holomorphy of $\evol_K$ (Lemma~III.4), combined with the complex 
linearity of the map 
$$ \Omega^1_h(M,\k) \to C^\infty([0,1],\k), \quad 
\alpha \mapsto \gamma^*\alpha  $$
implies that $\im(\delta)$ is a complex Lie group. 
As in Theorem~II.2(2), we see that 
$$ d \: {\cal O}_*(M,\k) \to T_0(\im(\delta)) $$
is an isomorphism of Lie algebras if ${\cal O}_*(M,\k)$ is endowed with the 
pointwise Lie bracket. 

{\bf Regularity:} 
It remains to verify the regularity of $G$, 
i.e., the holomorphy of the map 
$\evol_G$. As in the real case, we see that 
$$ \evol_G(\xi)(m) =  \evol_K(\xi^m) = \evol_G^\vee(\xi,m), $$
which is smooth by (Lemma~A.6(3)). Since $\evol_K$ is holomorphic 
(Lemma~III.4) and $\xi^m$ is complex linear in $\xi$ and holomorphic in $m$, 
it follows that $\evol_G^\vee \: C^\infty(I,\g) \times M \to K$ is holomorphic, 
which implies that $\evol_G$ is holomorphic. 
\qed

As in the real case, there is a natural situation where $\im(\delta)$ 
is a submanifold, namely if $M$ is one-dimensional and simply connected. 
In this case, each $\k$-valued holomorphic $2$-form on $M$ vanishes, 
so that all holomorphic $1$-forms satisfy the Maurer--Cartan
equation, and if $M$ is simply connected, Theorem~III.1 implies that 
$\im(\delta) = \Omega^1_{\rm dR,h}(M,\k)$, which is in particular a submanifold. 
If $M$ has no boundary, then the Riemann Mapping Theorem implies that 
it is isomorphic to $\C$, the unit disc 
$\Delta := \{ z \in \C \: |z| < 1\}$, or the Riemann sphere $\hat\C \cong \SS^2$. 

\Corollary III.9. For each regular complex Lie group $K$ and each $1$-connected 
complex curve $M$ without boundary, 
the group ${\cal O}_*(M,K)$ carries a regular complex 
Lie group structure for which 
$$ \delta \: {\cal O}_*(M,K) \to \Omega^1_h(M,\k) $$
is biholomorphic and ${\cal O}(M,K) \cong K \ltimes {\cal O}_*(M,K)$ 
carries a regular complex Lie group structure compatible with evaluations 
and the compact open topology. 
\qed

\Remark III.10. Let $M := \hat \C$ be the Riemann sphere. 
That the space $\Omega^1_h(\hat \C,\C)$ is trivial (which is well-known) 
can be seen as
follows. Each holomorphic $\C$-valued $1$-form $\alpha$ is closed, hence exact
because $\hat\C$ is simply connected. Now there exists a holomorphic
function $f \: \hat\C\to \C$ with $df = \alpha$, and since $f$ is
constant, we get $\alpha = 0$. {}From this and the Hahn--Banach Theorem, 
we directly get 
$\Omega^1(\hat\C,\k) = \{0\}$ for any locally convex space~$\k$. 

In view of Corollary~III.9 and the fact that $\delta$ is injective on 
${\cal O}_*(\hat \C,K)$, we now derive 
$${\cal O}_*(\hat \C,K) = \{\1\}$$ 
for any complex Lie group $K$, and this is independent of whether 
$K$ has non-constant holomorphic functions $K \to \C$ or not. 
In particular, we see that there is no non-constant holomorphic 
function from $\hat \C$ to any Lie group of the form 
$E/\Gamma_E$, where $\Gamma_E$ is a discrete subgroup of the complex locally 
convex space $E$. 
\qed

We extract the following version of the Regular Value Theorem from 
Gl\"ockner's Implicit Function Theorem ([Gl03]):  

\Theorem III.11. Let $M$ be a locally convex manifold, 
$N$ a Banach manifold, $F \: M \to N$ a smooth map and $n_0 \in N$. 
Assume that for each $m \in M$ with $F(m) = n_0$ there exists a continuous 
linear splitting of the tangent map 
$$ T_m(F) \: T_m(M) \to T_{n_0}(N). $$
Then $F^{-1}(n_0)$ is a split submanifold of $M$. 

If, in addition, $M$ and $N$ are complex manifolds and $F$ is holomorphic, 
then $F^{-1}(n_0)$ is a complex split submanifold of $M$. 

\Proof. Since the property of being a submanifold is local, it suffices to 
show that each $m_0 \in F^{-1}(n_0)$ has an open neighborhood $U$ for which 
$U \cap F^{-1}(n_0)$ is a submanifold of $U$. In particular, we may 
assume that $M$ is an open subset of a locally convex space 
$X\cong T_{m_0}(M)$. 
In view of the continuity of $F$, we may choose $U$ in such a way that 
$F(U)$ is contained in a chart domain in $N$, so that we may further assume 
that $N\cong T_{n_0}(N)$ is a Banach space. 

Fix a continuous linear splitting 
$\sigma \: N \to X$ of the tangent map $T_{m_0}(F)$. 
Then $Y := \ker T_{m_0}(F)$ is a closed subspace of $X$ and the map 
$$ Y \times N \to X, \quad (y,v) \mapsto y + \sigma(v) $$
is a linear topological isomorphism. 
We may therefore assume that $X = Y \times N$, write 
$m_0$ accordingly as $(y_0, e_0)$, and that $T_{m_0}(F) \: X = Y \times N \to N$
is the linear projection onto $N$. 

Now Theorem~2.3 in [Gl03] implies 
the existence of an open neighborhood 
$U$ of $m_0$ in $Y \times N$ and a diffeomorphism 
$\theta \: U \to \theta(U)$
onto some open neighborhood of $(y_0,n_0)$ 
in $Y \times N$ with 
$$ \theta(a,b) =  (a, F(a,b))\quad \hbox{ for } \quad (a,b) \in U. $$
This implies that 
$$ F^{-1}(n_0) \cap U = \theta^{-1}(Y \times \{n_0\}) $$
is a smooth submanifold of $U$. 
The remaining assertions are immediate from loc.~cit. 
\qed

The following theorem is the second main result of this section. 

\Theorem III.12. Let $M$ be a non-compact connected complex curve 
without boundary. Assume further that 
$\pi_1(M)$ is finitely generated and that $K$ is a complex Banach--Lie group. 
Then the group ${\cal O}_*(M,K)$ carries a regular complex 
Lie group structure for which 
$$ \delta \: {\cal O}_*(M,K) \to \Omega^1_h(M,\k) $$
is biholomorphic onto a complex submanifold, 
and ${\cal O}(M,K) \cong K \ltimes {\cal O}_*(M,K)$ carries a 
regular complex Lie group structure compatible with evaluations. 

\Proof. In view of Theorem~III.8, it suffices to show that 
$\im(\delta)$ is a complex submanifold. 

First we recall that the fundamental group 
$\pi_1(M)$ is free, because this is true for all non-compact surfaces without 
boundary. Let 
$$\gamma_1, \ldots, \gamma_r \: [0,1] \to M $$
be piecewise smooth loops in the base point $m_0$ such that 
$[\gamma_1],\ldots, [\gamma_r]$ are free generators of $\pi_1(M,m_0)$. 
Then the map 
$$ \Hom(\pi_1(M), K) \to K^r, \quad 
\chi \mapsto (\chi([\gamma_1]), \ldots, \chi([\gamma_r])) \leqno(3.1) $$
is a bijection. 

Since the Maurer--Cartan equation is trivially satisfied for holomorphic 
$1$-forms on a complex curve (cf.\ Corollary~III.9), we have 
$\im(\delta) = P^{-1}(\1)$
for the map 
$$ P \: \Omega^1_h(M,\k) \to K^r, \quad 
\alpha \mapsto (\per_\alpha^{m_0}([\gamma_1]), \ldots, \per_\alpha^{m_0}([\gamma_r])) $$
(Theorem~III.1). Since 
$$ \per_\alpha^{m_0}([\gamma]) = \evol_K(\gamma^*\alpha) $$
depends holomorphically on $\alpha$ for each piecewise smooth curve 
$\gamma$ on $M$ (Lemma~III.4), $P$ is a holomorphic map from the complex 
Fr\'echet space 
$\Omega^1_h(M,\k)$ to the complex Banach manifold~$K^r$. 

To see that $P^{-1}(\1)$ is a submanifold, we have to verify the assumptions of 
Theorem~III.11. 
The Behnke--Stein Theorem ([Fo77, Satz 28.6]) 
implies that each group homomorphism $\pi_1(M) \to \C$ 
can be realized by integration against a holomorphic $1$-form. 
Hence there exist holomorphic $1$-forms 
$\beta_1, \ldots, \beta_r \in \Omega^1_h(M,\C)$ with 
$$ \int_{\gamma_i} \beta_j = \delta_{ij}. \leqno(3.2) $$
We define a linear map  
$$ \sigma \: \k^r \to \Omega^1_h(M,\k), \quad 
(x_1,\ldots, x_r) \mapsto \sum_{j = 1}^r \beta_j \cdot x_j $$
whose continuity follows from Lemma~A.5(3). 

To verify for $\alpha \in P^{-1}(\1)$ 
that the map $T_\alpha(P)$ has a continuous linear section, we consider the map 
$$ \sigma_\alpha \: \k^r \to \Omega^1_h(M,\k), \quad 
x \mapsto \alpha + \Ad(f)^{-1}.\sigma(x) = \delta(f) + \Ad(f)^{-1}.\sigma(x), $$
where $f \in {\cal O}_*(M,K)$ is the unique function with $\delta(f)= \alpha$ 
(Theorem~III.1). 
As $f$ is fixed, $\sigma_\alpha$ is a continuous affine map, hence in 
particular holomorphic. From Remark~I.6(a) we further know that 
$$ \per_{\sigma_\alpha(x)}^{m_0} = \per_{\sigma(x)}^{m_0}, $$
so that $P\circ \sigma_\alpha = P \circ \sigma$. 

In view of $\per_\beta^{m_0}([\gamma]) = \evol_K(\gamma^*\beta)$, 
the differential of the map $\beta \mapsto \per_\beta^{m_0}$ in $0$ is given by 
$$ T_\1(\evol_K)(\gamma^*\beta) 
= \int_0^1 \gamma^*\beta = \int_\gamma \beta $$
(Lemma~A.5(1)). Therefore 
$$ T_0(P)(\beta) = \Big(\int_{\gamma_1} \beta, \ldots, \int_{\gamma_r} \beta\Big), $$
considered as an element of $\k^r$. From 
$$ \int_{\gamma_i} \sigma(x_1,\ldots,x_r) 
= \sum_{j = 1}^r \int_{\gamma_i} \beta_j \cdot x_j = x_i, $$
we derive 
$T_0(P) \circ \sigma = \id_{\k^r}.$
We further have 
$$ T_\alpha(P) \circ \Ad(f)^{-1} \circ \sigma 
= T_\alpha(P) \circ T_0(\sigma_\alpha) 
= T_0(P \circ \sigma_\alpha)
= T_0(P \circ \sigma)
= T_0(P)  \circ \sigma = \id_{\k^r}. $$
Hence $\Ad(f)^{-1} \circ \sigma$ is a continuous linear section of $T_\alpha(P)$. 
Since $\alpha \in P^{-1}(\1)$ was arbitrary, Theorem~III.11 implies that 
$P^{-1}(\alpha)$ is a complex 
submanifold of $\Omega^1_h(M,\k)$. 
\qed

\Corollary III.13. Let $\Sigma$ be a compact complex curve, $F \subeq \Sigma$ a finite 
set and $M := \Sigma \setminus F$. Then, for each Banach--Lie group $K$, 
the group  ${\cal O}(M, K)$ 
carries a regular complex 
Lie group structure compatible with the compact open topology 
and with evaluations. 

In particular, for each Banach--Lie group $K$, the topological group 
${\cal O}(\C^\times, K)$ carries a compatible Lie group structure. 

\Proof. To apply Theorem~III.12, it suffices to verify that $\pi_1(M)$ is finitely 
generated, but this follows from the fact that $M$ is homotopic  
to a compact surface with $|F|$ boundary circles. 
\qed

\Example III.14. Let $M = \C^\times$ and $K = \GL_n(\C)$. Then we associate 
to each holomorphic function $\xi \: \C^\times \to \gl_n(\C)$ the $1$-form 
$$ \alpha = \xi(z) dz. $$
Now $\delta(f) = \alpha$ is equivalent to the requirement that $f$ is a solution 
of the linear differential equation 
$$ f'(z) = f(z) \xi(z). $$

If, for example, $\xi(z) = z^{-1} A$ for a matrix $A$, then the differential 
equation reads 
$$ f'(z) = z^{-1} A f(z), \leqno(3.3) $$
and the corresponding $1$-form on the group $(\C^\times,\cdot)$ is invariant. 
A solution of (3.3) exists if and only if 
$$ f(z) = e^{\log z \cdot A} $$
is well-defined. The corresponding period homomorphism is 
$$ P(\alpha) \: \Z \to \GL_n(\C), \quad n \mapsto e^{2\pi i n A}. $$
Therefore $\alpha \in P^{-1}(\1)$ is equivalent to 
$e^{2\pi i A} = \1$, which is equivalent to the diagonalizability of 
$A$ and $\Spec(A) \subeq \Z$. 

Note that on the subspace $M_n(\C){dz \over z}$ the matrix $\1 \in \GL_n(\C)$ 
is not a regular value of $P$ because $\1$ is a not a regular value of 
the exponential function. Despite this fact, we have seen in the proof 
of Theorem~III.12 
that $\1$ is a regular value of the holomorphic function 
$$ P \: \Omega^1(\C^\times, M_n(\C)) \to \GL_n(\C), \quad \alpha \mapsto 
\per_\alpha. 
\qeddis 

\Remark III.15. Throughout the present section we considered only complex 
manifolds without boundary to make sure that the compact open topology 
on ${\cal O}(M,K)$ is the right one. If $M$ has non-empty boundary, 
all results remain true with respect to the finer smooth compact 
open topology. 
\qed

\sectionheadline{IV. Maps with values in special Lie groups} 

In this section we discuss the group $C^\infty(M,K)$ under the 
assumption that the universal covering group $\tilde K$ of $K$ 
is diffeomorphic to a locally convex space. This includes in particular 
all regular connected abelian Lie groups ([MT99], [GN07]), 
all finite-dimensional solvable Lie groups and many interesting 
projective limits of Lie groups ([HoNe06]). 

The starting point is the result that $C^\infty(M,\tilde K)$ always 
carries a Lie group structure compatible with evaluations. 
Then we study the passage from $K$ to $\tilde K$, 
which is encoded in an exact sequence 
$$ \1 \to \pi_1(K) \to 
C^\infty(M,\tilde K) \to 
C^\infty(M,K) \to \Hom(\pi_1(M),\pi_1(K)) \to \1 $$ 
which shows that the Lie group 
$C^\infty(M,\tilde K)/\pi_1(K)$ can be identified with a 
normal subgroup of $C^\infty(M,K)$, which eventually leads to a Lie 
group structure on the whole group. 
We further show that, if the group 
$\Hom(\pi_1(M),\pi_1(K))$, endowed with the topology of pointwise 
convergence, is discrete, then this Lie group structure 
is compatible with the smooth compact open topology. 
To shed some light on these subtleties, we briefly discuss 
the groups $\Hom(A,\Gamma)$ for an abelian group $A$ and a 
discrete subgroup $\Gamma$ of a locally convex space, and this discussion 
leads to necessary conditions for the Lie group structure on 
$C^\infty(M,K)$ to be compatible with the smooth compact open topology 
if $K$ is abelian or finite-dimensional. 

\Proposition IV.1. Let $M$ be a finite-dimensional manifold, 
$K$ a Lie group with Lie algebra $\k$, 
$\phi_K \: K \to E$ a diffeomorphism onto a locally convex space  
and $G := C^\infty(M,K)$. 
Then 
$$ \phi_G \: G \to C^\infty(M,E), \quad f \mapsto \phi_K \circ f $$
is a homeomorphism which defines a manifold structure on $G$, 
and this turns $G$ into a Lie group compatible with evaluations. 
If, in addition, $K$ is regular, then $G$ is also regular. 

\Proof. It follows directly from the functoriality of the topology on 
$C^\infty(M,N)$ in $N$ that $\phi_G$ is a homeomorphism. We consider 
$(\phi_G, G)$ as a smooth  atlas of the topological group $G$. 

To see that the group operations of $G$ are smooth, let 
$m_K \: K \times K \to K$ denote the multiplication of $K$ and 
$\eta_K$ its inversion map. Then 
$$ \phi_K \circ m_K \circ (\phi_K^{-1} \times \phi_K^{-1}) \: 
E \times E \to E $$
is smooth, so that Lemma~A.4 implies that the induced map 
$$ (\phi_K \circ m_K \circ (\phi_K^{-1} \times \phi_K^{-1}))_* \: 
C^\infty(M,E \times E) \to C^\infty(M,E) $$
is smooth, and this means that the multiplication in $G$ is smooth. 
With a similar argument, we see that the inversion is also smooth. 
Hence $G$ is a Lie group. 

To calculate the Lie algebra of this group, we observe that for 
$m \in M$, we have for the multiplication in the chart $(\phi_G,G)$ 
$$ \eqalign{ (f *_G g)(m) 
&:= \phi_G\Big(\phi_G^{-1}(f)\phi_G^{-1}(g)\Big)(m) 
= \phi_K\big(\phi_K^{-1}(f(m))\phi_K^{-1}(g(m))\big) \cr
&= f(m) *_K g(m)= f(m) + g(m) + b_\k(f(m), g(m)) + \cdots, \cr} $$
where the Lie bracket in $\k \cong T_\1(K) \cong E$ satisfies 
$$ [x,y] = b_\k(x,y) - b_\k(y,x) $$
(cf.\ [GN07]). Hence, we accordingly have 
$\big(b_\g(f,g)\big)(m) = b_\k(f(m),g(m)),$ and thus 
$$ [f,g](m) = b_\g(f,g)(m) - b_\g(g,f)(m) 
= b_\k(f(m),g(m))-b_\k(g(m),f(m)) = [f(m),g(m)].$$
Therefore $\g = C^\infty(M,\k)$, endowed 
with the pointwise defined Lie bracket, is the Lie algebra of~$G$. 

The compatibility of the Lie group structure with evaluations 
follows directly from  
Proposition~I.2 and the definition of the manifold structure on $G$. 

Now we assume that $K$ is regular. With Lemma~A.6(3), we obtain for each 
$\xi \in C^\infty(I,\g)$ a curve $\gamma \: I \to C^\infty(M,K)$ by 
$\gamma(t)(m) := \Evol_K(\xi^m)(t)$, 
defining a smooth map $I \times M \to K$, hence a smooth curve in 
$G$ (Lemma~A.2) because 
$$ C^\infty(I,G) \cong 
C^\infty(I,C^\infty(M,E)) \cong 
C^\infty(I \times M,E) \cong 
C^\infty(I \times M,K). $$
Further, $\delta(\gamma^m) = \xi^m$ implies that 
the evolution map of $G$ is given by 
$\evol_G(\xi)(m) := \evol_K(\xi^m)$.  
Now the smoothness of $\evol_G$ follows from Lemma~A.6(3) 
and Lemma~A.2. 
\qed

\Theorem IV.2. Let $M$ be a finite-dimensional 
connected manifold, 
$K$ a connected Lie group with Lie algebra $\k$ whose universal 
covering group $\tilde K$ is diffeomorphic to a locally convex space. 
 Then the following assertions hold: 
\litem{(1)} $G := C^\infty(M,K)$ carries the structure of a Lie group 
compatible with evaluations. 
If $K$ is regular, then $G$ is regular. 
\litem{(2)} On the identity component $G_0$, 
this Lie group structure is compatible with the smooth compact open topology. 
In particular, the Lie group structure on $G$ is 
compatible with the smooth compact open topology if and only if 
$G_0$ is open with respect to the smooth compact open topology. 
\litem{(3)} If $M$ is $\sigma$-compact, then 
$\pi_0(G) \cong \Hom(\pi_1(M), \pi_1(K))$ with respect to the 
Lie group structure and the smooth compact open topology. 
\litem{(4)} If $\pi_1(M)$ is finitely generated or, more generally, 
$\Hom(\pi_1(M), \pi_1(K))$ 
is discrete with respect to the topology of pointwise convergence, then 
$G_0$ is also open in the 
smooth compact open topology, so that the Lie group structure on $G$ is 
compatible with this topology. 

\Proof. (1) First we apply Proposition~IV.1 to the Lie group 
$\tilde G_0 := C^\infty(M,\tilde K)$ to obtain a Lie group structure 
compatible with evaluations. 
Let $q_K \: \tilde K \to K$ denote the universal covering group. 
Since $\ker q_K$ is central in $\tilde K$, the conjugation action 
$$C_{\tilde K} \: \tilde K \times \tilde K \to \tilde K, 
\quad (x,y) \mapsto xyx^{-1} $$
factors through a smooth action 
$$\tilde C_K \: K \times \tilde K \to \tilde K, \quad (x,y) \mapsto xyx^{-1}. $$
We claim that the corresponding action of 
$G = C^\infty(M,K)$ on $\tilde G_0 = C^\infty(M,\tilde K)$ by 
$(f.g)(m) := \tilde C_K(f(m))(g(m))$, is an action by  
smooth automorphisms. To see this, first observe that if $q_M \: \tilde M \to M$ 
is the universal covering manifold of $M$, then 
$C^\infty(\tilde M,\tilde K)$ also carries a smooth Lie group 
structure for which we may identify $C^\infty(M,\tilde K)$ with a closed submanifold 
(corresponding to a closed vector subspace under the chart in Proposition~IV.1). 
Since each smooth map $f \: M \to K$ can be lifted to a smooth map 
$\tilde f \: \tilde M \to \tilde K$, the corresponding automorphism 
of $\tilde G_0$ coincides with the restriction of the 
automorphism $c_{\tilde f}$ of $C^\infty(\tilde M,\tilde K)$ to 
$C^\infty(M,\tilde K)$, hence is smooth. 

Since $\ker q_K \cong \pi_1(K) \subeq \tilde K$ is a discrete central subgroup 
of $\tilde K$ and therefore also of the Lie group $\tilde G_0$, the quotient 
$G_0 := \tilde G_0/\pi_1(K)$ carries a unique Lie group structure for which 
the quotient map $\tilde G_0 \to G_0$ is a covering (cf.\ [GN07]). 
This quotient map corresponds to the homomorphism 
$$ q_K^M \:  C^\infty(M,\tilde K)\to C^\infty(M,K), \quad 
f \mapsto q_K \circ f $$
which is equivariant with respect to the aforementioned action of the 
group $C^\infty(M,K)$ on $C^\infty(M,\tilde K)$ by conjugation. Hence 
$G_0 \cong \im(q_K^M)$ is a normal subgroup which carries a Lie group 
structure and the other group elements act by smooth automorphisms. 
This implies that the Lie group structure extends uniquely to all of 
$G$ in such a way that $G_0$ is the open identity component of $G$ 
(cf.\ [GN07]). We thus obtain a Lie group structure on $G$ for which 
$q_K^M \: \tilde G_0 \to G_0$ is the universal covering map, which 
implies that its Lie algebra also is $\g$. 
The evaluation map $\ev \: G \times M \to K$ can be obtained by factorization of 
the evaluation map $\tilde\ev \: \tilde G_0 \times M \to \tilde K$ because 
$$ q_K \circ \tilde\ev = \ev \circ (q_K^M \times \id_M), $$
which shows that it is smooth on $G_0 \times M$ and hence on all of 
$G \times M$ because it is multiplicative in the first argument. 

If, in addition, $K$ is regular, then Proposition IV.1 implies that 
$\tilde G_0$ is regular, and this easily implies that $G_0$ and $G$ 
are regular Lie groups (cf.\ [GN07]).

(2) The construction in Proposition~IV.1 implies that on 
$\tilde G_0 = C^\infty(M,\tilde K)$ the Lie group structure 
is compatible with the smooth compact open topology. 
Writing this group as a semidirect product 
$\tilde G_0 \cong C^\infty_*(M,\tilde K) \rtimes \tilde K$, we see that 
$G_0 \cong C^\infty_*(M,\tilde K) \rtimes K$, so that it suffices to see that 
the injective continuous homomorphism 
$$ q_K^M \: C^\infty_*(M,\tilde K) \to C^\infty_*(M,K) $$
of topological groups is a topological embedding with respect to the 
smooth compact open topology. 

To verify this claim, let $m \in M$ and $C$ be a relatively compact open 
neighborhood of $m$. To see that 
$$ (q_K^M)^{-1} \: q_K^M(C_*^\infty(M,\tilde K)) \to C^\infty_*(M,\tilde K) $$
is continuous in the smooth compact open topology, we note that 
for each $m \in \N$, the map 
$$ T^m(q_K) \: T^m\tilde K \to T^mK $$ 
is the universal covering morphism of the Lie group $T^mK$. If 
we can show that for all these coverings, the corresponding map 
$$ (q_{T^mK}^M)^{-1} \: q_{T^mK}^M(C_*(M,T^m \tilde K)) \to C_*(M,T^m\tilde K) $$
is continuous in the compact open topology, the corresponding 
assertion follows. Hence it suffices to show that 
$$ (q_K^M)^{-1} \: q_K^M(C_*(M,\tilde K))_c \to C_*(M,\tilde K)_c \leqno(4.1) $$
is continuous with respect to the compact open topology. 

Evaluation in the base point $m_0$ maps the 
subgroup $\pi_1(K) \cong \ker q_K$ to a discrete subgroup of $\tilde K$. 
Hence it is discrete as a subgroup of $C(M,\tilde K)_c$, so that for any 
sufficiently small $\1$-neighborhood $V$ in $C(M,\tilde K)$ we 
have $C_*(M,\tilde K) \cap V\pi_1(K) \subeq V$. 
Let $C$ be a compact connected subset of $M$ and 
$U$ be an open $\1$-neighborhood in $\tilde K$ with 
$$ W(C,U) \pi_1(K) \cap C_*(M,K) \subeq W(C,U) 
\quad \hbox{ and } \quad UU^{-1} \cap \pi_1(K) = \{\1\}. $$
We claim that for any $\tilde f \in C_*(M,\tilde K)$ 
the relation $q_K^M(\tilde f) \in W(C,q_K(U))$ implies 
$\tilde f \in W(C,U)$, and this implies the continuity of (4.1). 
Any such $\tilde f$ maps the connected set $C$ into the open subset 
$U\pi_1(K) = q_K^{-1}(q_K(U))$, and our assumption on $U$ implies that 
the open sets $Uz$, $z \in \pi_1(K)$, are pairwise disjoint. Hence there exists 
a $z \in \pi_1(K)$ with 
$\tilde f \in W(C,Uz) = W(C,U)z$, so that $\tilde f \in C_*(M,\tilde K)$ 
leads to $z = \1$, i.e., $\tilde f \in W(C,U)$. 

This completes the proof of (2). We also note that if $K$ is regular, 
then the assertion follows directly from Proposition~II.1, because both maps 
$$ \delta_1 \: C^\infty_*(M,\tilde K) \to \Omega^1(M,\k) 
\quad \hbox{ and } \quad  \delta_2 \: C^\infty_*(M, K) \to \Omega^1(M,\k) $$
are topological embeddings with 
$\delta_2 \circ q_K^M = \delta_1$. 

(3) The range of $q_K^M$ consists of all smooth maps $f \: M \to K$ lifting to 
maps $\tilde f \: M \to \tilde K$. If $m_0 \in M$ is a base point, this 
condition is equivalent to the condition that the homomorphism 
$$ \Gamma(f) \: \pi_1(M,m_0) \to [\SS^1, K] \cong \pi_1(K), \quad 
[\alpha] \mapsto [f\circ \alpha] $$
is trivial. In view of 
$$ [(f_1 \cdot f_2) \circ \alpha] 
= [(f_1 \circ \alpha) \cdot (f_2 \circ \alpha)] 
= [f_1 \circ \alpha]\cdot [f_2 \circ \alpha], $$
$\Gamma$ is a group homomorphism, and we obtain an  exact sequence 
$$ \1 \to \pi_1(K) \into 
C^\infty(M,\tilde K) \smapright{q_K^M} 
C^\infty(M,K) \sssmapright{\Gamma} \Hom(\pi_1(M),\pi_1(K)) \leqno(4.2) $$
of groups. 

We also note that $\Gamma$ is continuous with respect to the 
compact open topology on $C^\infty(M,K)$ and the topology of pointwise 
convergence on the abelian group $\Hom(\pi_1(M), \pi_1(K))$, 
which turns it into a totally disconnected group because $\pi_1(K)$ is 
discrete. This implies that $\im(q_K^M)$ also coincides with the 
arc-component of $\1$ with respect to the smooth compact open topology. 

To see that $\Gamma$ is surjective, let $\gamma \in \Hom(\pi_1(M),\pi_1(K))$,  
and consider the corresponding $\tilde K$-principal bundle 
$P_\gamma \:= \tilde M \times_\gamma \tilde K \to M$ over $M$. 
Since $\tilde K$ is contractible and $M$ is $\sigma$-compact and connected, 
hence paracompact, 
this bundle is topologically 
trivial, and hence also smoothly trivial ([MW06]), so that there exists 
a smooth function $\tilde f \: \tilde M \to \tilde K$ with 
$\tilde f(d.x) = \gamma(d) \tilde f(x)$ for 
$d \in \pi_1(M)$ and $x \in \tilde M$. Then $\tilde f$ induces a smooth 
function $f \: M \to K$ with $\Gamma(f) = \gamma$. 
As $\ker \Gamma$ is the arc-component of the identity for both topologies 
on $C^\infty(M,K)$, its surjectivity implies (3). 

(4) If $\pi_1(M)$ is finitely generated, the group $\Hom(\pi_1(M),\pi_1(K))$ 
is discrete with respect to the topology of pointwise convergence. 
This in turn implies that $G_0$ is open with respect to the smooth compact open 
topology. Hence (2) implies that the Lie group structure on $G$ is compatible 
with this topology because this is the case on $G_0$. 
\qed

\Theorem IV.3. Let $M$ be a finite-dimensional 
connected complex manifold, 
$K$ a connected complex Lie group with Lie algebra $\k$ whose universal 
covering group $\tilde K$ is diffeomorphic to a locally convex space. 
 Then the following assertions hold: 
\litem{(1)} $G := {\cal O}(M,K)$ carries a Lie group structure compatible 
with evaluations. If $K$ is regular, then $G$ is regular. 
\litem{(2)} On the identity component $G_0$, 
this Lie group structure is compatible with the smooth compact open topology. 
\litem{(3)} If $\pi_1(M)$ is finitely generated or $\Hom(\pi_1(M), \pi_1(K))$ 
is discrete with respect to the topology of pointwise convergence, then 
$G_0$ is also open in the 
smooth compact open topology, so that the Lie group structure on $G$ is 
compatible with this topology. 
\litem{(4)} If $M$ is Stein (hence $\sigma$-compact) and $K$ is Banach, then 
$\pi_0(G) \cong \Hom(\pi_1(M), \pi_1(K))$ holds with respect to the 
Lie group structure and the smooth compact open topology. 

\Proof. Using Lemma~A.4(3), (1)--(3) are proved exactly as in the real case, 
where the evaluation map is holomorphic by Proposition~I.2(2). 
We omit the details. 

(4) Here the crucial step is the surjectivity of the map $\Gamma$, for which 
we need that any bundle $P_\chi = M \times_\chi \tilde K \to K$ is 
holomorphically trivial. In the proof of Theorem~IV.2, we have argued that it is 
topologically trivial. Since $M$ is assumed to be Stein and $K$ Banach, 
the Oka principle (cf.\ [Rae77, Th.~2.1]) asserts that this 
bundle has a holomorphic section. The remaining arguments are similar to the real case. 
\qed

\Remark IV.4. (a) It is quite plausible that under the assumptions of 
Theorem~IV.2, the abelian group $\pi_1(K)$ is torsion free 
because each torsion element would lead to a fixed point free action 
of a cyclic group on the locally convex space $\tilde K$. 

If $K$ is finite-dimensional, such an action does not exist (cf.\ [Sm41]), 
and if $K$ is regular and abelian, we know anyhow that $\pi_1(K)$ is a subgroup of the 
additive group of $\k$, hence torsion free. 

(b) Theorem~IV.2 applies in particular to all 
finite-dimensional Lie groups $K$ which are diffeomorphic to some vector space. 
These groups are isomorphic to semidirect products of the form 
$$ R \rtimes \tilde\SL_2(\R)^m, $$
where $R$ is a simply connected solvable Lie group and 
$\tilde\SL_2(\R)$ is the universal covering group of $\SL_2(\R)$, which 
is diffeomorphic to $\R^3$. 

Interesting infinite-dimensional Lie groups diffeomorphic to locally convex 
spaces can also be found among the simply connected pro-solvable pro-Lie groups 
(see [HoNe06] for more details). 

(c) Theorem~IV.3 applies to all 
finite-dimensional complex Lie groups $K$ which are biholomorphic to a vector space, 
which is equivalent to $K$ being solvable. 
\qed

\subsection{Maps with values in abelian Lie groups} 

We have seen in Theorem~IV.2  that the smooth compact open topology on 
$C^\infty(M,K)$ is compatible with the Lie group topology if and only if 
the arc-component of the identity is open. To get a better understanding 
of what is going wrong if $\pi_1(M)$ is not finitely generated, we now take a closer 
look at the relevant facts on abelian groups. 

Let $K$ be a regular connected abelian Lie
group, hence of the form $K \cong \k/\Gamma_K$, where $\k$ is a Mackey complete 
locally convex space and 
$\Gamma_K \subeq \k$ is a discrete subgroup isomorphic to $\pi_1(K)$ 
(cf.\ [MT99], [GN07]). We write $q_K \: \k \to K$ for the quotient map with kernel 
$\Gamma_K$. Then $\tilde K \cong \k$ is a locally convex space and 
Theorem~IV.2 applies. 

Since $\k$ is abelian, $\alpha \in \Omega^1(M,\k)$ satisfies the 
Maurer--Cartan equation if and only if $\alpha$ is closed and for 
each $m_0 \in M$ we have 
$$\per^{m_0}_\alpha \: \pi_1(M,m_0) \to K, \quad 
 [\gamma] \mapsto q_K\Big(\int_\gamma \alpha\Big). $$
Therefore a  closed $1$-form $\alpha$ is integrable (this is called 
logarithmically exact in [Pa61]) if and only if 
all its periods are contained in $\Gamma_K$. Let 
$$ Z^1_{\rm dR}(M,\k,\Gamma_K) := \im(\delta) \subeq Z^1_{\rm dR}(M,\k) $$
denote the subgroup of all $1$-forms satisfying this condition (Theorem~I.5). 
Proposition~II.1 now implies that 
$$ \delta \: C^\infty_*(M,K) \to Z^1_{\rm dR}(M,\k,\Gamma_K) \leqno(4.3) $$ 
is an isomorphism of topological groups, inducing an 
isomorphism of the path components of the identity  
$$ C^\infty_*(M,K)_a \cong C^\infty_*(M,\k) 
\to dC^\infty(M,\k) $$
because we have 
$$\Gamma \: C^\infty(M,K) \to \Hom(\pi_1(M), \Gamma_K), \quad 
\Gamma(f)([\gamma]) = \int_\gamma \delta(f). $$
We conclude that the path component of the identity in 
$C^\infty(M,K)$ is open if and only if 
$$ H^1_{\rm dR}(M,\k,\Gamma_K) := Z^1_{\rm dR}(M,\k,\Gamma_K)/dC^\infty(M,\k)$$
is a discrete subgroup of the topological vector space $H^1_{\rm dR}(M,\k) 
= Z^1_{\rm dR}(M,\k)/dC^\infty(M,\k)$. 
In view of the de Rham Theorem (cf.\ [KM97]) and the Hurewicz Homomorphism, 
we have isomorphisms of groups 
$$ H^1_{\rm dR}(M,\k) \cong H^1_{\rm sing}(M,\k) \cong \Hom(H_1(M),\k) \cong
\Hom(\pi_1(M),\k) $$
because $H_1(M) \cong \pi_1(M)/[\pi_1(M),\pi_1(M)]$
(cf.\ [Bre93]). By restriction, we thus obtain the isomorphism of groups  
$$ H^1_{\rm dR}(M,\k,\Gamma_K) \cong \Hom(\pi_1(M),\Gamma_K), \leqno(4.4) $$
and, in view of Theorem~IV.2, it is interesting to see when this actually 
is an isomorphism of topological groups, resp., when the 
subgroup $H^1_{\rm dR}(M,\k,\Gamma_K)$ is discrete. 

\Lemma IV.5. If $M$ is connected and $\sigma$-compact and 
$\k$ is a Fr\'echet space, then the natural map 
$$ \Phi \: H^1_{\rm dR}(M,\k) \to \Hom(H_1(M),\k) $$
is an isomorphism of Fr\'echet spaces, where the space 
$\Hom(H_1(M),\k)$ carries the topology of pointwise convergence. 

\Proof. First we observe that it is a continuous 
bijection of Fr\'echet spaces, where the topology of pointwise convergence
on $\Hom(H_1(M),\k)$ defines a Fr\'echet space structure because 
the group $H_1(M)$ is countably generated (cf.\ [Ne04, Prop.~IV.9]). 
Now the Open Mapping Theorem ([Ru73]) implies that 
$\Phi$ is an isomorphism of Fr\'echet spaces. 
\qed

In view of the preceding lemma, we are left with the question when 
$\Hom(A,\Gamma_K)$ is discrete in $\Hom(A,\k)$ for an abelian group $A$. 
The following theorem provides some information on the structure of $\Gamma_K$. 

\Theorem IV.6. {\rm(Sidney)} Countable 
discrete subgroups of locally convex spaces are free. 

\Proof. If $\Gamma$ is a discrete subgroup of the locally convex space $\k$, 
then there exists a continuous seminorm $p$ on $\k$ such that 
$\inf \{ p(\gamma) \: 0 \not=\gamma \in \Gamma \} > 0.$
If $\k_p$ is the completion of the normed space $\k/p^{-1}(0)$, it
follows that $\Gamma$ embeds as a subgroup of $\k_p$ whose intersection
with a sufficiently small ball is trivial, and therefore $\Gamma$ is
realized as a discrete subgroup of some Banach space. 
This implies that every discrete subgroup 
of a locally convex space is isomorphic to a discrete subgroup of 
some Banach space. Now the assertion follows from 
Sidney's Theorem that countable discrete
subgroups of Banach spaces are free ([Si77, p.983]).  
\qed

\Lemma IV.7. Let $A$ be a countable abelian group. Then the following
are equivalent: 
\litem{(1)} $\Hom(A,\Gamma)$ is discrete in $\Hom(A,E)$ for each
discrete subgroup $\Gamma$ of a locally convex space~$E$. 
\litem{(2)} $\Hom(A,\Gamma)$ is discrete in $\Hom(A,E)$ for one non-zero 
discrete subgroup $\Gamma$ of some locally convex space~$E$. 
\litem{(3)} $\Hom(A,\Z)$ is discrete in $\Hom(A,\R)$. 
\litem{(4)} $\Hom(A,\Z)$ is finitely generated. 

\Proof. Let $A_1 \subeq A$ denote the intersection of all kernels of
homomorphisms $A \to \Z$. Then we have an embedding 
$$ A/A_1 \into \Z^{\Hom(A,\Z)}, \quad a + A_1 \mapsto (\chi \mapsto \chi(a)) $$
and the group $A/A_1$ is countable. As all countable subgroups of
groups of the form $\Z^I$ are free ([Fu70, Th.\ 19.2]), 
we have 
$$  A \cong A_1 \oplus A_2, $$
where $\Hom(A_1,\Z) = \0$ and $A_2 \cong \Z^{(J)}$ is free 
(cf.\ [Fu70, Cor.~19.3] for this result due to K.~Stein). 
In particular, we have $\Hom(A,\Z) \cong \Hom(A_2, \Z) \cong \Z^J$, 
and this group is finitely generated if and only if $J$ is finite. 

If $\chi \: A \to \Gamma$ is a homomorphism, then 
$\chi(A)$ is a countable discrete subgroup of the locally convex space
$\k$, hence free by Sidney's Theorem IV.6. Therefore the homomorphisms 
$\chi(A) \to \Z$ separate points, which implies that 
$A_1 \subeq \ker \chi$. We conclude that 
$$ \Hom(A,\Gamma) \cong \Hom(A_2,\Gamma) \cong \Hom(\Z^{(J)},\Gamma) 
\cong \Gamma^J. $$
The topology of pointwise convergence on $\Hom(A,\Gamma)$
corresponds to the product topology on $\Gamma^J$. Hence this group is
discrete if and only if either $\Gamma = \{0\}$ or $J$ is finite. 
In particular, $\Hom(A,\Z)$ is discrete if and only if $J$ is
finite. 
\qed

{}From the proof of Lemma IV.7 we obtain in particular that 
$\Hom(A,\Z) \cong \Z^J$, where $J$ is a countable set. 
If $J$ is finite, then $\Z^J \cong \Z^{|J|}$ is discrete, and otherwise it is 
isomorphic to $\Z^\N$, which is not discrete. 

\Theorem IV.8. Let $M$ be a connected $\sigma$-compact 
manifold, $\k$ a Fr\'echet space, 
$\Gamma_K \subeq \k$ a non-trivial discrete subgroup, and 
$K := \k/\Gamma_K$. Then the following are equivalent: 
\litem{(1)} The Lie group structure on $C^\infty(M,K)$ from {\rm Theorem~IV.2} is 
compatible with the smooth compact open topology. 
\litem{(2)} The arc-component $C^\infty(M,K)_a = q_K^M(C^\infty(M,\k))$ 
is open with respect to the smooth compact open topology. 
\litem{(3)} $H^1_{\rm dR}(M,\k,\Gamma_K)$ is a discrete subgroup of
$H^1_{\rm dR}(M,\k)$. 
\litem{(4)} $H^1(M,\Z)$ is finitely generated. 

\Proof. The equivalence of (1) and (2) follows from Theorem~IV.2(2). 

The equivalence of (1) and (3) follows from the discussion preceding Lemma~IV.5. 
To prove the equivalence between (3) and (4), we recall 
from (4.4) and  Lemma~IV.5 the isomorphism
$$ H^1_{\rm dR}(M,\k,\Gamma_K) \cong \Hom(H_1(M),\Gamma_K) $$
of topological groups. 
As $H_1(M)$ is a countable abelian group (cf.\ [Ne04, Prop.~IV.8]), 
the equivalence of (3) and
(4) follows from $\Gamma_K \not=\{0\}$, 
Lemma IV.7, and $\Hom(H_1(M),\Z) \cong H^1(M,\Z)$. 
\qed

\subheadline{The complex case} 

Now we assume that $M$ is a complex manifold and $K \cong \k/\Gamma_K$ 
is a regular abelian complex Lie group, 
i.e., $\k$ is a complex Mackey complete space. 

\Theorem IV.9. Let $M$ be a connected $\sigma$-compact 
complex manifold without boundary, $\k$ a complex Fr\'echet space,  
$\Gamma_K \subeq \k$ a non-zero discrete subgroup and $K := \k/\Gamma_K$.
 Then 
$$ \delta \: {\cal O}_*(M,K) \to Z^1_{\rm dR,h}(M,\k,\Gamma_K) $$
is an isomorphism of topological groups, and the following are
equivalent: 
\litem{(1)} The Lie group structure on ${\cal O}(M,K)$ is compatible 
with the compact open  topology. 
\litem{(2)} The arc-component ${\cal O}_*(M,K)_a = q_K^M({\cal O}(M,\k))$ 
is open with respect to the compact open topology. 
\litem{(3)} $H^1_{\rm dR,h}(M,\k,\Gamma_K) 
= \{[\alpha] \: \alpha \in \MC_h(M,\k), 
(\forall \gamma \in C^\infty(\SS^1,M)) \int_\gamma \alpha \in \Gamma_K\}$ 
is a discrete subgroup of $H^1_{\rm dR,h}(M,\k)$. 
\nin  If, in addition, $M$ is a Stein manifold, then {\rm(1)-(3)} are 
equivalent to: 
\litem{(4)} $H^1(M,\Z)$ is finitely generated. 

\Proof. The equivalence of (1)--(3) is shown precisely as in 
the real case. Suppose, in addition, that $M$ is a Stein manifold.
Then the Oka principle ([Gr58], Satz I, p.~268) implies that each continuous map 
$M \to \C^\times$ is homotopic to a holomorphic map. 
As the homotopy classes $[M,\C^\times]$ are classified by 
$$\Hom(\pi_1(M),\pi_1(\C^\times)) \cong \Hom(\pi_1(M),\Z) 
\cong H^1(M,\Z), $$ 
each homomorphism $\pi_1(M) \to \Z$ arises as 
$\pi_1(f)$ for a holomorphic map 
$f \: M \to \C^\times$. Hence each class in 
$H^1(M,\Z)$ is represented by a 
holomorphic $1$-form. In view of Lemma~IV.5, this implies that 
$$  H^1_{\rm dR,h}(M,\C,\Z) 
= H^1_{\rm dR}(M,\C,\Z) \cong \Hom(H_1(M),\Z) $$
as topological groups, 
and the discreteness of this group is equivalent to (4) (Proposition~IV.8). 
\qed

\Remark IV.10. The assumption that $M$ is Stein is crucial in (4) above, 
because, in general, not every homomorphism $H_1(M) \to \C$ is 
represented by integration of a holomorphic $1$-form. 
A typical example is given by a complex torus 
$M = \C/D$, where $\Z^2 \cong D \subeq \C$ is a discrete subgroup. 
Since holomorphic functions on $M$ are constant, 
each holomorphic $1$-form on $M$ is a constant multiple of $dz$. 
Therefore $H^1_{\rm dR,h}(M,\C) = \C[dz]$ is one-dimensional, and 
the homomorphism 
$\pi_1(M) \cong D \to \C$ corresponding to $\lambda \cdot dz$ 
is given by $d \mapsto \lambda d$. The group homomorphism 
$D \to \C, d \mapsto \oline d$ is not represented by 
integration of a holomorphic $1$-form. 
\qed

\Remark IV.11. If $M$ is a non-compact Riemann surface, then $M$ is Stein 
and $\pi_1(M)$ is a free group, so that 
$H_1(M)$ is a free abelian group. Therefore Theorem~IV.9  
implies that for any abelian complex Fr\'echet--Lie group 
of the form $K=\k/\Gamma_K$, the group ${\cal O}(M,K)$ is a Lie group 
with respect to the compact open topology if and only if 
the free abelian group $H_1(M)$ has finite rank. 
\qed

\Remark IV.12. (a) If $M$ is a compact K\"ahler manifold and $\k$ abelian, 
then $\Omega^1_h(M,\k) =  \MC_h(M,\k)$
because each holomorphic $1$-form is automatically closed (cf.\ [We80]). 

(b) For any compact complex manifold we have ${\cal O}_*(M,\k) = \{0\}$ because 
all holomorphic functions are constant. Therefore the Lie group structure on 
${\cal O}_*(M,K)$ is discrete. 
\qed

\subsection{Maps with values in finite-dimensional Lie groups} 

\Remark IV.13. (a) Let $K$ be a connected finite-dimensional Lie group 
whose universal covering group $\tilde K$ is diffeomorphic to a vector space, 
which is equivalent to $\k/\rad(\k) \cong \sL_2(\R)^m$ for some $m \in \N_0$ 
and this in turn is equivalent to the maximal 
compact subgroup $T \subeq K$ being a torus (cf.\ [HoNe06]). 
Let $d := \dim T$. 
Since $T$ is a maximal compact subgroup, 
$\pi_1(K) \cong \pi_1(T) \cong \Z^d$ is a free 
abelian group (cf.\ [Ho65]). Therefore 
$$ \Hom(\pi_1(M), \pi_1(K)) \cong \Hom(H_1(M), \pi_1(T)) 
\cong H^1(M,\Z)^d. $$
If $K$ is not simply connected, then $d > 0$, so that 
this group is discrete if and only if 
$H^1(M,\Z)$ is finitely generated (Lemma~IV.7). 

If this is the case, then Theorem~IV.2 implies that 
 the Lie group structure on $C^\infty(M,K)$ is 
compatible with the smooth compact open topology. 
If $H^1(M,\Z)$ is not finitely generated, then 
$$C^\infty_*(M,\tilde T) = C^\infty_*(M,\tilde K) \cap C^\infty_*(M,T)$$ 
is not an open subgroup of $C^\infty_*(M,T) \subeq C^\infty_*(M,K)$, and therefore 
$C^\infty_*(M,\tilde K)$ is not an open subgroup of $C^\infty_*(M,K)$. 
Now the topological decomposition 
$$ C^\infty(M,K) \cong C^\infty_*(M,K) \rtimes K $$
implies that the arc-component of the identity in 
$C^\infty(M,K)$ is not open, hence that the Lie group structure 
is not compatible with the smooth compact open topology. 

(b) If $K$ is a finite-dimensional complex Lie group, the Levi decomposition 
implies that the condition that 
$\tilde K$ is diffeomorphic to a vector space is equivalent to $K$ being solvable. 

Suppose that $H^1(M,\Z)$ is not finitely generated
 and let $T \subeq K$ be a real maximal torus. 
Then the inclusion $T \into K$ extends to a holomorphic Lie group morphism 
$T_\C \to K$. Note that for $T \cong \T^n$ the universal complexification is 
$T_\C \cong (\C^\times)^n$.
If the map $T_\C \into K$ is an embedding, then we also have an embedding 
$$ T_\C \into {\cal O}(M,T_\C) \into {\cal O}(M,K). $$
If $M$ is a Stein manifold, 
then ${\cal O}(M,T_\C)$ is not a Lie group 
because ${\cal O}(M,\tilde T_\C) \cong {\cal O}(M,T_\C)_a$ 
is not open (Theorem~IV.9). As in (a) above, 
this implies that ${\cal O}(M,\tilde K)$ is not open in ${\cal O}(M,K)$, 
and therefore that ${\cal O}(M,K)$ is not a Lie group with respect to the 
compact open topology. 
\qed

\sectionheadline{V. Some strange properties of the exponential map} 

In this subsection, we collect some interesting properties of the exponential 
function of the groups ${\cal O}(M,K)$ on a finite-dimensional 
complex manifold $M$ to a complex Banach--Lie group $K$, which is simply given by 
$\exp(\xi) := \exp_K \circ \xi$, where $\exp_K$ is the exponential function 
of $K$. 

\Proposition V.1. Let 
$M$ be a complex manifold which has non-constant holomorphic 
functions, and $K$ be a complex connected Banach-Lie group.
 If $K \not= \exp_K \k$, then the image of the exponential function 
of ${\cal O}(M,K)$ is not an identity neighborhood. 

\Proof. {\bf Step 1:} First we claim the existence of a holomorphic function 
$\ell \: M \to \C$ with real part unbounded from above. Suppose that such a function does 
not exist. Replacing $f$ by $if$, $-f$ and $-if$, we conclude that 
for each holomorphic function $f \: M \to \C$ the functions 
$\Re f$ and $\Im f$ are bounded, and hence that $f$ is bounded. 
If $f$ is non-constant, then $f(M)$ is an open subset of $\C$, hence has a boundary 
point $z_0 \not\in f(M)$. But then the function $(f - z_0)^{-1}$ is unbounded, 
a contradiction. We therefore find a holomorphic function $\ell \: M \to \C$ and a sequence 
$x_n \in M$ with $\Re \ell(x_n) \to \infty$. 

{\bf Step 2:} Let 
$$ K_1 := \{ f(1) \: f \in {\cal O}(\C,K), f(0) = \1\}. $$
Then $K_1$ is a subgroup of $K$, because it is the homomorphic image of the 
subgroup ${\cal O}_*(\C,K)$ under the evaluation map in $1$. 
If $k = \exp_K x$ for some $x \in \k$, then the map 
$f(z) := \exp_K(zx)$ satisfies $f(0) = \1$ and $f(1) = k$. Hence $K_1 \supeq 
\exp_K \k$, 
and since the connected Banach--Lie group $K$ is generated by $\exp_K \k$, 
we obtain $K = K_1$. 

{\bf Step 3:} Let $k \in K \setminus \exp_K \k$. In view of the 
preceding paragraph, there exists 
a holomorphic map $f \: \C \to K$ with $f(1) = k$ and $f(0) = \1$. 
We define $h_n(x) := f(e^{\ell(x)-\ell(x_n)})$. Then 
$h_n(x_n) = f(1) \not\in \exp_K \k$, so that $h_n$ is not contained in the 
image of the exponential function of ${\cal O}(M,K)$. 
On the other hand $h_n \to \1$ uniformly on compact subsets of $M$, hence in 
${\cal O}(M,K)$. 
\qed

\Corollary V.2. Let $M$ be a complex manifold with non-constant holomorphic 
functions and $K_1\leq K$ a Banach--Lie subgroup whose exponential 
function is not surjective. 
Then there exist $0$-neighborhoods in ${\cal O}(M,\k)$ 
whose image under the exponential function is not an identity neighborhood 
in ${\cal O}(M,K)$. 

\Proof. Let $U_\k \subeq \k$ be an open $0$-neighborhood 
which is relatively compact and for which 
$\exp_K\res_{2 U_\k}$ is a diffeomorphism onto its open image, satisfying 
$$ \exp_K(U_\k) \cap K_1 = \exp_K(U_\k \cap \k_1). \leqno(5.1) $$
Pick $m_0 \in M$ and a compact neighborhood $C$ of $m_0$. 
Then we consider the identity neighborhood 
$W(C,\exp_K(U_\k))$ of ${\cal O}(M,K)$. 
Let $\k_1 \leq \k$ be the Lie algebra of $K_1$ 
and observe that $W(C, U_\k)$ is a $0$-neighborhood in 
${\cal O}(M,\k)$. 

In view Proposition~V.1, each identity neighborhood 
in ${\cal O}(M,K_1)$ 
contains a holomorphic function $h \: M \to K_1$, 
not contained in the image of the exponential function of ${\cal O}(M,K_1)$. 
Suppose that $h = \exp \xi=\exp_K \circ \xi$ holds for some holomorphic 
function $\xi\:M\to\k$, contained in $W(C,U_\k)$. 
Then the injectivity of $\exp_K\res_{U_\k}$ and (5.1) imply  that 
$\xi(C) \subeq U_\k \cap \k_1$. Since $f$ is holomorphic, we obtain 
$\xi(M) \subeq \k_1$, contradicting the construction of $h$. 
Therefore $\exp(W(C,U_\k))$ is not an identity neighborhood 
in ${\cal O}(M,K)$. 
\qed

The preceding corollary implies in particular that the exponential function of \break 
${\cal O}(M,\GL_n(\C))$ is not locally surjective for any Stein manifold $M$. 
The following lemma shows that it is locally injective. 

\Lemma V.3. If $M$ is a connected complex manifold, then the exponential 
function $$\exp \: {\cal O}(M,\k) \to {\cal O}(M,K)$$ is locally injective. 

\Proof. Let $C \subeq M$ be a non-empty compact subset, 
$U_\k \subeq \k$ be an open $0$-neighborhood on which the exponential function 
$\exp_K \: \k \to K$ is injective, and define 
$$ U := W(C,U_\k) = \{ \xi \in {\cal O}(M,\k) \: (\forall x \in C)\ \xi(x) \in U_\k\}. $$
Then $U$ is an open $0$-neighborhood in ${\cal O}(M,\k)$. 
If $\xi,\eta \in U$ satisfy $\exp \xi = \exp \eta$, then  
the injectivity of $\exp_K$ on $U_\k$ implies that 
$\xi\res_C = \eta\res_C$, and since $M$ is connected, we obtain 
$\xi =\eta$ by analytic continuation. 
\qed

\sectionheadline{Appendix. Technical tools} 

\Lemma A.1. Let $M, N$ 
and $L$ be locally convex manifolds, $f \in C^\infty(M \times N,L)$ 
and put $f_x(y) := f(x,y)$. 
Then the map 
$f^\vee \: M \to C^\infty(N,L), x \mapsto f_x$ 
is continuous. 

\Proof. {\bf Step 1:} (cf.\ [Ne01, Lemma~III.2]) For Hausdorff spaces 
$M$, $N$ and $L$ and $f \in C(M \times N,L)$, the map 
$f^\vee \: M \to C(N,L)_c$, 
is continuous: 
Suppose that $f_x \in W(K,U)$ for some compact subset $K \subeq N$ and 
some open subset $U \subeq L$, i.e.,
$\{x \} \times K \subeq f^{-1}(U)$. Since $f^{-1}(U)$ is an open
subset of $M \times N$ and $\{x\} \times K \subeq M \times N$ is
compact, there exists an open neighborhood $O \subeq M$ of $x$ such
that $O \times K \subeq f^{-1}(U)$. This means that 
$x \in O \subeq \{ p \in M \: f_p \in W(K,U)\}$, which proves the
assertion. 

\nin{\bf Step 2:} $f^\vee$ is continuous. For each $k \in \N$ we have 
a natural product decomposition \break $T^k(M \times N) \cong T^k(M) \times T^k(N)$, 
so that Step 1 implies the continuity of the maps 
$$ M \to C(T^k(N),T^k(L))_c, \quad x \mapsto T^k(f_x) $$
for each $k \in \N$. In view of the definition of the topology on 
$C^\infty(N,L)$, this proves that $f^\vee$ is continuous. 
\qed

\Lemma A.2. Let $N$ and $M$ be smooth manifolds. 
\litem{(1)} If $E$ is a locally convex space and $f \in C^\infty(N \times M,E)$, 
then $f^\vee \: N \to C^\infty(M,E)$ is smooth. 
\litem{(2)} If $M$ is compact (possibly with boundary) 
and $K$ is a Lie group, then for each smooth map 
$f \in C^\infty(N \times M,K)$, the map 
$f^\vee \: N \to C^\infty(M,K)$ is smooth with respect to the 
natural Lie group structure on $C^\infty(M,K)$. 

\Proof. (1) We may w.l.o.g.\ assume that 
$N$ is an open convex subset of a locally convex space $X$ and identify 
$T(N)$ with $N \times X$. First we show that $f^\vee$ is $C^1$ with tangent map 
$$\Psi \: T(N) \to C^\infty(M,T(E)),\quad 
\Psi_m(v)(n) = T_{(m,n)}(f)v $$
whose continuity follows from Lemma~A.1 and 
the smoothness of the tangent map 
$$T(f) \: T(N) \times T(M) \to T(E).$$ 

Fix $(x,h) \in T(N)$. For a sufficiently small $\eps > 0$ the map
$$ ]-\eps, \eps[ \times [0,1] \to C^\infty(M,E), \quad 
(t,u) \mapsto \Psi(x + uth,h) $$
is continuous by Lemma~A.1. Therefore 
$$ ]-\eps, \eps[ \to C^\infty(M, E), \quad 
t \mapsto \int_0^1 \Psi(x + uth,h)\ du $$
is continuous, and so 
$$ \lim_{t \to 0} {1\over t} \big( f^\vee(x + th) - f^\vee(x)\big) 
= \lim_{t \to 0} \int_0^1 \Psi(x + uth,h)\ du
= \int_0^1 \Psi(x,h)\ du = \Psi(x,h). $$
Thus $T(f^\vee)(x,h) = \Psi(x,h)$, and the continuity of $\Psi$ implies
that $f^\vee$ is $C^1$. 

Applying this argument to the tangent map $T(f^\vee)$, we see that 
$T(f^\vee)$ is also $C^1$, so that $f^\vee$ is $C^2$. Proceeding inductively, it 
follows that $f^\vee$ is smooth. 

(2) To see that $f^\vee$ is smooth in a neighborhood of some 
$n_0 \in N$, it suffices to prove the smoothness of the map 
$n \mapsto f^\vee(n_0)^{-1} f^\vee(n)$, so that we may assume that 
$f^\vee(n_0) = \1$. 

Let $(\phi_K,U_K)$ be a $\k$-chart of $K$ and 
$(\phi,U)$ the corresponding $C^\infty(M,\k)$-chart of the group 
$C^\infty(M,K)$, given by 
$$ U := C^\infty(M,U_K) \quad \hbox{ and } \quad \phi(\xi) = \phi_K \circ \xi $$
(cf.\ Theorem~I.3). Then the continuity of $f^\vee$ implies that 
$f^\vee$ maps a neighborhood $U_N \subeq N$ of $n_0$ into $U$, we may assume that 
$f(N \times M) \subeq U$, and we have to 
show that $\phi \circ f^\vee\res_{U_N}$ is smooth. 
Since 
$$ (\phi \circ f^\vee)(n)(m) = \phi_K(f(n,m)) = (\phi_K \circ f)^\vee(n)(m), $$
we have 
$\phi \circ f^\vee = (\phi_K \circ f)^\vee \: N \to C^\infty(M,\k),$
and the smoothness of this map follows from (1). 
\qed

\Lemma A.3. Let $\K \in \{\R,\C\}$, $N$ be a locally convex 
smooth $\K$-manifold,  
$M$ a finite-dimensional smooth 
$\K$-manifold (without boundary in case $\K = \C$)
and $E$ a topological $\K$-vector space. 
Then the following assertions hold: 
\litem{(1)} For $\K = \R$, a map $f \: N \to C^\infty(M,E)$ is smooth 
if and only if the map 
$$f^\wedge \: N \times M \to E, \quad f^\wedge(n)(m) := f(n,m) $$
is smooth. The map $\Psi \: C^\infty(N,C^\infty(M,E)) \to 
C^\infty(N \times M, E), 
f \mapsto f^\wedge$ is an isomorphism of topological vector spaces. 
\litem{(2)} For $\K = \C$, a map $f \: N \to {\cal O}(M,E)$ is holomorphic  
if and only if $f^\wedge$ is holomorphic. The 
map $\Psi \: {\cal O}(N,{\cal O}(M,E)) \to 
{\cal O}(N \times M, E), 
f \mapsto f^\wedge$ is an isomorphism of topological vector spaces. 

\Proof. [Gl04, Prop.~12.2] directly implies (1). 
To verify (2), we first observe that 
the Cauchy Formula implies that on the closed subspace 
${\cal O}(M,E)$, uniform convergence on compact subsets 
implies uniform convergence 
of all partial derivatives on compact subsets. Hence the inclusion map 
$$ {\cal O}(M,E) \into C^\infty(M,E) $$
is continuous and therefore a topological embedding. 
In this sense the compact open topology on ${\cal O}(M,E)$ coincides 
with the $\C$-smooth compact open topology, which is used in 
[Gl04]. Therefore (2) follows from [Gl04, Prop.~12.2] for $\K= \C$ or 
by observing that the map $\Psi$ in (1) maps the closed subspace 
${\cal O}(N,{\cal O}(M,E))$ of $C^\infty(N,C^\infty(M,E))$ 
homeomorphically onto ${\cal O}(N \times M, E)$. 
\qed

\Lemma A.4. Let $E_1$ and $E_2$ be locally convex spaces, $U_j \subeq E_j$ 
open subsets, and 
$\phi \: U_1 \to U_2$ be a smooth map. 
\litem{(1)} The map 
$$ \phi_* \: C^\infty(M,U_1) \to C^\infty(M,U_2), 
\quad f \mapsto \phi \circ f $$
is continuous. 
\litem{(2)} If, in addition, $M$ is compact or $U_j = E_j$ for $j =1,2,$ so that 
the subsets $C^\infty(M,U_j)$ are open in $C^\infty(M,E_j)$, then the map 
$\phi_*$ is smooth. 
\litem{(3)} If, in addition to the assumptions in {\rm(2)}, 
$E_1$ and $E_2$ are complex and $\phi$ is holomorphic, then 
$\phi_*$ is holomorphic. 

\Proof. (1) The continuity of $\phi_*$ follows directly from the definition 
of the topology and the continuity of left compositions with respect 
to the compact open topology. 

(2) Assume that $M$ is compact or $U_j = E_j$ for $j =1,2$. 
In view of Lemma~A.2, the smoothness of $\phi_*$ 
follows from the smoothness of the corresponding map 
$$ C^\infty(M,U_1) \times M \to U_2, \quad (f,m) \mapsto \phi(f(m)) = 
\phi \circ \ev(f,m), $$
where $\ev \: C^\infty(M,U_1) \times M \to U_1$ is the 
smooth evaluation map 
(Proposition~I.2). 

(3) In view of (2), it remains to show that the differentials of $\phi_*$ 
are complex linear, which follows from 
$$ d(\phi_*)(f)(\xi)(x) = d\phi(f(x))\xi(x). 
\qeddis 

 The following two lemmas collect some technical smoothness properties
of regular Lie groups. 

\Lemma A.5. For a connected finite-dimensional 
smooth manifold $M$ (with boundary) and a regular Lie group $K$ with Lie algebra~$\k$, 
the following assertions hold: 
\litem{(1)} The map 
$\Evol_K \: C^\infty([0,1],\k) \to C^\infty_*([0,1],K)$
is a diffeomorphism with 
$$ T_0(\Evol_K)(\xi)(t) = \int_0^t \xi(s)\, ds 
\quad \hbox{ and } \quad 
T_0(\evol_K)(\xi) = \int_0^1 \xi(s)\, ds. $$
Its inverse is 
$$ \delta \: C^\infty_*([0,1],K) \to C^\infty([0,1],\k) \quad 
\hbox{ with } \quad  T_\1(\delta)(\xi) = \xi'. $$
\litem{(2)} The action of $K$ on 
$\Omega^1(M,\k)$ by $\Ad(k).\alpha := \Ad(k) \circ \alpha$ 
is smooth. 
\litem{(3)} The multiplication map $\Omega^1(M,\R) \times \k \to \Omega^1(M,\k), 
(\alpha,x) \mapsto \alpha \cdot x$ is continuous. 

\Proof. (1) First we show that $\Evol_K$ is smooth. 
For $\gamma \in C^\infty([0,1],K)$ and $\gamma_s(t) := \gamma(st)$ we have 
$\delta(\gamma_s)(t) = s\delta(\gamma)(st) = S(s,\delta(\gamma))(t),$
where 
$$ S \: [0,1] \times C^\infty([0,1],\k) \to C^\infty([0,1],\k),
 \quad S(s,\xi)(t) = s\xi(st) $$
is smooth by Proposition~I.2 and Lemma~A.2. 
For $\delta(\gamma) = \xi$ we have 
$$ \Evol_K(\xi)(s) = \gamma(s) = \gamma_s(1) = \evol_K(S(s,\xi)),$$
showing that the map 
$C^\infty([0,1],\k) \times [0,1] \to K, (\xi,s) \mapsto \Evol_K(\xi)(s)$
is smooth, and this implies that $\Evol_K$ is smooth 
(Lemma~A.2). 

To see that $\delta$ is smooth, we write 
$\delta(\gamma)(t) = \kappa_K(\gamma'(t))$. Since $\kappa_K$ is smooth, 
the assertion follows from the smoothness of the homomorphism of 
Lie groups  
$$ C^\infty([0,1],K) \to C^\infty([0,1],TK), \quad 
\gamma \mapsto \gamma', $$
Lemma A.2 and the smoothness of the evaluation map 
of $C^\infty([0,1],TK)$ (Theorem~I.3). 

>From $\delta \circ \Evol_K = \id_{C^\infty([0,1],\k)}$ and the smoothness 
of $\Evol_K$, we derive that 
$$ T_\1(\delta) \circ T_0(\Evol_K) = \id_{C^\infty([0,1],\k)}. $$
Using the Chain Rule, we obtain directly 
$T_\1(\delta)(\xi)(t) = \xi'(t),$
and since $T_\1(\delta)$ is injective on $C^\infty_*([0,1],\k)$, the tangent space of 
$C^\infty_*([0,1],K)$ in the constant function $\1$, we get 
$$T_0(\Evol_K)(\xi)(t) = \int_0^t \xi(s)\, ds. $$
Now $\ev_1 \circ \Evol_K = \evol_K$ leads to the asserted formula 
for $T_0(\evol_K)$. 

(2) Since $\Omega^1(M,\k)$ is a closed subspace of 
$C^\infty(TM,\k)$, it suffices to observe that the action of 
$K$ on $C^\infty(TM,\k)$, given by $k.f := \Ad(k) \circ f$ is smooth. 
In view of Lemma~A.2,  it suffices to show that the map 
$$K \times C^\infty(TM,\k) \times TM \to \k, \quad 
(k,f,m) \mapsto \Ad(k).f(m) $$
is smooth, 
which in turn follows from the smoothness of the adjoint action of $K$ on $\k$ and 
the smoothness of the evaluation map of $C^\infty(TM,\k)$ (Proposition~I.2). 

(3) With the same argument as in (2), it 
suffices to show that 
$$ C^\infty(TM,\R) \times \k \to C^\infty(TM,\k), \quad (f,x) \mapsto f\cdot x $$
is smooth. This in turn follows from the smoothness of the map 
$$ C^\infty(TM,\R) \times TM \times \k \to \k, \quad (f,v,x) \mapsto f(v)\cdot x 
= \ev(f,v)\cdot x $$ 
(Proposition~I.2, Lemma~A.2). 
\qed

\Lemma A.6. Let $M$ be a connected finite-dimensional 
smooth manifold (with boundary) and $K$ a regular Lie group with Lie algebra~$\k$. 
\litem{(1)} If $\gamma \: [0,1] \to M$ is a piecewise smooth curve, then 
the map 
$$ \Omega^1(M,\k) \mapsto K, \quad \alpha \mapsto \evol_K(\gamma^*\alpha) $$
is smooth. 
\litem{(2)} Let $(\phi,U)$ be a chart of $M$ for which 
$\phi(U)$ is a convex $0$-neighborhood 
and $\gamma_x(t) := \phi^{-1}(t\phi(x))$. Then the map 
$$ \Omega^1(M,\k) \times U \to K, \quad 
(\alpha, x) \mapsto \evol_K(\gamma_x^*\alpha) $$
is smooth. 
\litem{(3)} For $\xi \in C^\infty(I \times M,\k)$ 
put $\xi^m(t) := \xi(t,m)$. Then the map 
$$ \gamma \: I \times M \to K, \quad (t,m) \mapsto \Evol_K(\xi^m)(t) $$
is smooth with 
$$ \delta(\gamma)_t(m) := \gamma(t,m)^{-1}.{d\over dt} \gamma(t,m) = \xi^m(t), $$
and  the map 
$$ \evol_G^\wedge \: C^\infty(I,\g) \times M \to K, \quad 
(\xi, m) \mapsto \evol_K(\xi^m) $$
is also smooth. 

\Proof. (1) This follows 
from the smoothness of $\evol_K$ and the fact that for each 
smooth path $\eta \: [a,b] \to M$ the map 
$$ \Omega^1(M,\k) \to C^\infty([a,b],\k), \quad \alpha \mapsto \eta^*\alpha 
= \alpha \circ T\eta $$
is continuous and linear, hence smooth. 

(2) Since $K$ is regular, we have to show that the map 
$$ \Omega^1(M,\k) \times U \to C^\infty([0,1],\k), \quad 
(\alpha, x) \mapsto \gamma_x^*\alpha $$
is smooth. 
In view of Lemma~A.2, this 
follows from the smoothness of the map 
$$ \Omega^1(U,\k) \times U \times [0,1] \to \k, \quad 
(\alpha, x,t) \mapsto (\gamma_x^*\alpha)_t = \alpha_{\gamma_x(t)}\gamma_x'(t),$$
which is a consequence of the smoothness of the evaluation map of 
$C^\infty(TU,\k)$ (Proposition~I.2) and of the map 
$U \times [0,1] \to TM,(x,t) \mapsto \gamma_x'(t)$. 

(3) First we recall from Lemma~A.3 
that for $\g := C^\infty(M,\k)$ we have 
$$ C^\infty(I,\g) 
\cong C^\infty(I \times M,\k) 
\cong C^\infty(M, C^\infty(I,\k)) $$
as topological vector spaces. 
In this sense, we consider each $\xi \in C^\infty(I,\g)$ as a smooth map 
$I \times M \to \k$. In particular, $\xi^m \in C^\infty(I,\k)$, 
$\evol_K(\xi^m) \in K$, and the map 
$\xi^\vee \: M \to C^\infty(I,\k), m \mapsto \xi^m$ is smooth.  
Hence the smoothness of $\gamma$ follows from 
$$ \gamma(t,m) =\Evol_K(\xi^m)(t) 
= \ev\circ (\Evol_K(\xi^\vee(m)),t)
= \ev\circ ((\Evol_K \circ \xi^\vee) \times \id_I)(m,t) $$
because the evaluation map of $C^\infty(I,K)$ is smooth (Theorem~I.3). 
The formula for $\delta(\gamma)$ follows immediately from the definition. 

To see that $\evol_G^\wedge$ is smooth, we first recall that 
$\evol_K$ is smooth. Hence it suffices to observe that the map 
$$ C^\infty(I \times M, \k) \times M \to C^\infty(I,\k), \quad 
(\xi,m) \mapsto \xi^m $$ 
is smooth because it corresponds to the evaluation map of the space 
$C^\infty(M,C^\infty(I,\k))$ (cf.\ Proposition~I.2). 
\qed

\def\entries{

\[Beg87 Beggs, E. J., {\it The de Rham complex on infinite dimensional manifolds}, 
Quart. J. Math. Oxford (2) {\bf 38} (1987), 131--154 

\[Bre93 Bredon, G.\ E., ``Topology and Geometry,'' Graduate Texts in
Mathematics {\bf 139}, Springer-Verlag, Berlin, 1993 

\[Fo77 Forster, O., ``Riemannsche Fl\"achen,'' Heidelberger Taschenb\"ucher, 
Springer-Verlag, 1977 

\[Fu70 Fuchs, L., ``Infinite Abelian Groups, I,'' 
Acad. Press, New York, 1970 

\[Gl02a Gl\"ockner, H., {\it Lie group structures on quotient groups
and universal complexifications for infinite-dimensional Lie groups},
J. Funct. Anal.  {\bf 194}  (2002),  no. 2, 347--409

\[Gl02b ---, {\it Infinite-dimensional Lie groups without completeness 
restrictions}, in ``Geometry and Analysis on Finite-
and Infinite-Dimensional Lie Groups,'' A.~Strasburger et al Eds., 
Banach Center Publications {\bf 55}, Warszawa 2002; 53--59 

\[Gl03 ---, {\it Implicit Functions from Topological Vector Spaces to Banach 
Spaces}, Israel Journal Math. {\bf 155} (2006), 205--252 

\[Gl04 ---, {\it Lie groups over non-discrete topological fields},
arXiv:math.GR/0408008 

\[Gl05 ---, {\it H\"older continuous homomorphisms between infinite-dimensional 
Lie groups are smooth}, J. Funct. Anal. {\bf 228:2} (2005),  419-444

\[GN07 Gl\"ockner, H., and K.-H. Neeb, ``Infinite-dimensional Lie Groups,'' 
book in preparation 


\[Gr58 Grauert, H., {\it Analytische Faserungen \"uber holomorph-vollst\"andigen R\"aumen}, 
Math. Ann. {\bf 135} (1958), 263--273

\[GR77 Grauert, H., and R. Remmert, {\it Theorie der Steinschen R\"aume}, 
Grundlehren der math. Wissenschaften {\bf 227}, Springer, 1977

\[Gro68 Grothendieck, A., {\it Classes de Chern et representations lineaires}, 
in ``Dix Expos\'es sur la Cohomologie des Sch\'emas,'' North Holland, 
Amsterdam; Masson, Paris, 1968 

\[Ha82 Hamilton, R., {\it The inverse function theorem of Nash and
  Moser}, Bull. Amer. Math. Soc. {\bf 7} (1982), 65--222

\[Ho65 Hochschild, G., ``The Structure of Lie Groups,'' Holden Day, San 
Francisco, 1965 

\[HoNe06 Hofmann, K.~H., and K.-H.\ Neeb, 
{\it Pro-Lie groups which are infinite-dimensional Lie groups}, submitted 

\[KT68 Kamber, F., and Ph.\ Tondeur, ``Flat manifolds,'' Lecture Notes Math. 
{\bf 67}, Springer-Verlag, 1968 

\[KM97 Kriegl, A., and P.\ Michor, ``The Convenient Setting of
Global Analysis,'' Math.\ Surveys and Monographs {\bf 53}, Amer.\
Math.\ Soc., 1997 

\[Mi80 Michor, P. W., ``Manifolds of Differentiable
Mappings,'' Shiva Publishing, Or\-pington, Kent (U.K.), 1980 

\[MT99 Michor, P., and J.~Teichmann, {\it Description of infinite
dimensional abelian regular Lie groups}, J. Lie Theory {\bf 9:2}
(1999), 487--489 

\[Mil58 Milnor, J., {\it On the existence of a connection with curvature zero}, Comment. Math. 
Helv. {\bf 32} (1958), 215--223 

\[Mil84 ---,  {\it Remarks on infinite-dimensional Lie groups},
pp.~1007--1057; In:  DeWitt, B., Stora, R. (eds), 
``Relativit\'{e}, groupes et topologie II (Les Houches, 1983), North Holland, Amsterdam, 1984 

\[MW06 M\"uller, Chr., and Chr. Wockel, {\it Equivalences of smooth and continuous 
principal bundles with infinite-dimensional structure group}, 
math.DG/0604142 

\[Ne01 Neeb, K.-H., {\it Representations of infinite dimensional
  groups}, pp.~131--178; in ``Infinite Dimensional K\"ahler Manifolds,'' 
Eds. A. Huckleberry, T. Wurzbacher, DMV-Seminar {\bf 31}, 
Birkh\"auser Verlag, 2001 

\[Ne04 ---, {\it Current groups for non-compact manifolds and their
central extensions}, pp.\ 109--183; in ``Infinite Dimensional Groups and Manifolds''
IRMA Lect. Math. Theor. Phys. {\bf 5}, de Gruyter Berlin 2004

\[Ne06 ---, {\it Towards a Lie theory of locally convex 
groups}, Jap. J. Math. 3rd ser. {\bf 1:2} (2006), 291--468 

\[Pa61 Palais, R., {\it Logarithmically exact differential forms}, Proc. Amer. 
Math. Soc. {\bf 12} (1961), 50--52 

\[Rae77 Raeburn, I., {\it The relationship between a commutative Banach
algebra and its maximal ideal space}, J. Funct. Anal. {\it 25} (1977),
366--390 

\[Ram65 Ramspott, K. J., {\it Stetige und holomorphe Schnitte in
B\"undeln mit homogener Faser}, Math. Z. {\bf 89} (1965), 234--246 

\[Ru73 Rudin, W., ``Functional Analysis,'' McGraw Hill, 1973

\[Sch03 Schlichenmaier, M., {\it Higher genus affine algebras of 
Krichever-Novikov type}, Mosc. Math. J.  {\bf 3:4} (2003),  1395--1427

\[Si77 Sidney, S. J., {\it Weakly dense subgroups of Banach spaces},
Indiana Univ. Math. Journal {\bf 26:6} (1977), 981--986 

\[Sm41 Smith, P. A., {\it Fixed-point theorems for periodic transformations}, 
Amer. J. Math. {\bf 63} (1941), 1--8

\[We80 Wells, R.~O., Jr., ``Differential Analysis on Complex Manifolds,'' 
Second edition, Graduate Texts in Mathematics {\bf 65}, Springer-Verlag, New York-Berlin, 1980

\[Wo05 Wockel, Chr., {\it The Topology of Gauge Groups}, 
submitted; math-ph/0504076 

\[Wo06 ---, {\it Smooth Extensions and Spaces of Smooth and Holomorphic 
Mappings}, J.~Geom. Symm. Phys. {\bf 5} (2006), 118--126 

}

\references
\dlastpage

\bye